\documentclass[11pt]{amsart}
\usepackage{amsmath,amsfonts}

\newtheorem{thm}{Theorem}[subsection]
\newtheorem{lem}[thm]{Lemma}
\newtheorem{prop}[thm]{Proposition}
\newtheorem{dfn}[thm]{Definition}

\newtheorem{rmk}[thm]{Remark}

\newtheorem{ex}[thm]{Example}

\newtheorem*{ack}{Acknowledgments}
\newtheorem*{fothm}{Fukaya-Oh Theorem}

\def\del{\partial}

\def\R{\mathbb{R}}
\def\C{\mathbb{C}}
\def\Z{\mathbb{Z}}

\def\on{\operatorname}
\def\CD{\mathcal D}
\def\CF{\mathcal F}
\def\CG{\mathcal G}
\def\CL{\mathcal L}
\def\CN{\mathcal N}
\def\CO{\mathcal O}
\def\CS{\mathcal S}
\def\CT{\mathcal T}

\def\ol{\overline}
\def\intHom{{\mathcal Hom}}
\def\dist{{d}}

\def\CA{\mathcal A}
\def\CB{\mathcal B}
\def\CJ{\mathcal J}
\def\CY{\mathcal Y}

\def\modules{mod}
\def\hom{{hom}}
\def\ob{{Ob}}
\def\ch{{Ch}}
\def\tw{{Tw}}

\def\CC{\mathcal C}

\def\supp{\on{supp}}
\def\fU{\mathfrak U}
\def\CM{\mathcal M}

\def\vareps{{\varepsilon}}
\def\CE{\mathcal E}

\def\fY{\mathfrak Y}
\def\wt{\widetilde}

\def\codim{\on{codim}}

\begin{document}

\title{Constructible Sheaves and the Fukaya Category}
\author{David Nadler and Eric Zaslow}
\begin{abstract} 
Let $X$ be a compact real analytic manifold, and let $T^*X$ be its cotangent bundle.
Let $Sh(X)$ be the triangulated dg category 
of bounded, constructible complexes of sheaves on $X$.
In this paper, we develop a Fukaya $A_\infty$-category $Fuk(T^*X)$ 
whose objects are exact,
not necessarily compact Lagrangian
branes in the cotangent bundle. 
We write $Tw Fuk(T^*X)$ for the $A_\infty$-triangulated envelope of $Fuk(T^*X)$
consisting of twisted complexes of Lagrangian branes.
Our main result is that $Sh(X)$ quasi-embeds into
$Tw Fuk(T^*X)$ as an $A_\infty$-category. 
Taking cohomology gives an
embedding of the corresponding derived categories.
\end{abstract}

\maketitle


{\tiny
\tableofcontents
}

\section{Introduction}
In this paper, we study the relationship between two natural invariants
of a real analytic manifold $X$. 
The first is the Fukaya category of Lagrangian submanifolds of the cotangent bundle $T^*X$.
The second is the derived category of constructible sheaves on $X$ itself.
The two are naively related by the theory of linear differential equations --
that is, the study of modules over the ring $\CD_X$
of differential operators on $X$.
On the one hand, Lagrangian cycles in $T^*X$ play a prominent role in 
the microlocal theory of $\CD_X$-modules. 
On the other hand, in the complex setting, 
the Riemann-Hilbert correspondence identifies 
regular, holonomic $\CD_X$-modules with constructible sheaves.
In what follows, we give a very brief account of what we mean by the Fukaya category of $T^*X$
and the constructible derived category of $X$, and then state our main result.

Roughly speaking, 
the Fukaya category of a symplectic manifold
is a category whose
objects are Lagrangian submanifolds 
and whose morphisms and compositions
are built out of the quantum intersection theory of Lagrangians.
This
is encoded by the moduli space of pseudoholomorphic maps from polygons
with prescribed Lagrangian boundary conditions.
Since $T^*X$ is noncompact, there are many choices to be made
as to which Lagrangians to allow and how to obtain well-behaved moduli spaces
of pseudoholomorphic maps. One approach is to insist that the Lagrangians
are compact. With this assumption, the theory is no more difficult than
that of a compact symplectic manifold. One perturbs the Lagrangians
so that their intersections are transverse, and then convexity arguments guarantee 
compact moduli spaces.

Our version of the Fukaya category $Fuk(T^*X)$ includes 
both compact and noncompact exact Lagrangians.
We work with
exact Lagrangians that have well-defined limits at infinity.
To make this precise, we consider a compactification of $T^*X$,
and assume that the closures of our Lagrangians are subanalytic
subsets of the compactification.
Two crucial geometric statements follow from this assumption. First,
the boundaries of our Lagrangians
are Legendrian
subvarieties of the divisor at infinity. Second, for any metric on the fibers of $T^*X$,
its restriction
to one of our Lagrangians
has no critical points near infinity.
These facts allow us to make sense of ``intersections at infinity"
by restricting our perturbations to those which are normalized geodesic flow
near infinity
for carefully prescribed times.
Given suitable further perturbations (which are available
in intended applications),
 we then obtain compact
moduli spaces of pseudoholomorphic maps.
The resulting Fukaya category $Fuk(T^*X)$ has many of the usual properties
that one expects from both a topological and categorical perspective.

The second invariant of the real analytic manifold $X$ which we consider
is the derived category $D_c(X)$ of constructible sheaves on $X$ itself.
This is a triangulated category which encodes the topology of subanalytic subsets
of $X$. To give a sense of the size of $D_c(X)$, its
Grothendieck group is the group of constructible functions on $X$ -- that is, functions
which are constant along some subanalytic stratification,
for example a triangulation.
Examples of objects of $D_c(X)$ include closed submanifolds equipped
with flat vector bundles. More generally, we have the so-called
standard and costandard objects associated to 
a locally closed submanifold $Y\subset X$ equipped with a flat vector bundle $\CE$.
Informally, one may think of the standard object
as 
the complex of singular cochains on $Y$ with values in $\CE$,
and the costandard
as the complex of relative singular cochains on $(Y,\del Y)$ with values in $\CE.$
A key observation is that morphisms between these objects are naturally
the singular cohomology of certain subsets of $X$
with values in flat vector bundles.

One formulation of our main result is the following.
As we outline below, it may be viewed as a categorification of the characteristic cycle
construction.

\begin{thm}
Let $X$ be a real analytic manifold. Then there is a canonical embedding 
of derived categories
$$
D_c(X)\hookrightarrow DFuk(T^*X).
$$
\end{thm}

The result reflects an underlying  quasi-embedding of dg and $A_\infty$-categories.
Further arguments shows that this is in fact a quasi-equivalence~\cite{Nad}.

The remainder of the introduction is divided into several parts.
In the section immediately following, we discuss motivations for our main result result from the long-developing theory of 
microlocal geometry.
In the section after that, we explain the general outline
of the proof of our main result.
Finally, we speculate on possible applications
in the context of mirror symmetry.


\subsection{Microlocal geometry}

The main result of this paper has a natural place in the context of microlocal geometry.
Broadly speaking, sheaf theory on a real analytic manifold $X$ 
may be thought of as a tool to understand local 
analytic and topological
phenomena
and how they assemble into global phenomena. Many
aspects of the theory are best understood from a microlocal perspective, or in other words
as local phenomona on the cotangent bundle $T^*X$.
We collect here a short account of some results from this subject
that naturally point toward our main result. What we present is not intended
to be an exhaustive overview of the subject. For that we refer the reader
to the book
of Kashiwara-Schapira~\cite{KS}.
It contains many original results,
presents a detailed development of the subject,
and includes historical notes and a comprehensive bibliography.

Our main result may be viewed as a categorification of the characteristic cycle
construction for real constructible sheaves introduced by Kashiwara~\cite{Kindex}.
(For foundational material on microlocal constructions
such as the singular support, see Kashiwara-Schapira~\cite{KSmicrolocal}.)
Given a constructible complex of sheaves $\CF$ on $X$,
its characteristic cycle $CC(\CF)$ is a conical Lagrangian
cycle in $T^*X$ (with values in the pullback of the orientation sheaf of $X$)
which encodes the singularities of the original complex. The multiplicity
of $CC(\CF)$ at a given covector is the Euler characteristic of the local Morse groups of
the complex with respect to the covector. If a covector is not in the support
of $CC(\CF)$, it means that there is no obstruction to propagating local sections
of $\CF$ in the direction of the covector. So for example, the characteristic cycle
of a flat vector bundle on $X$ is the zero section in $T^*X$
with multiplicity the dimension of the vector bundle. More generally,
the characteristic
cycle of a flat vector bundle on a closed submanifold is the conormal bundle to the submanifold
with multiplicity the dimension of the vector bundle.

As mentioned earlier, the Grothendieck group of the constructible derived category $D_c(X)$
is the space of constructible functions on $X$.
The characteristic cycle construction descends to an isomorphism from constructible functions
to the group of conical Lagrangian cycles in $T^*X$. From this vantage point, there are many
results that might lead one to our main result. First, there is the index formula of 
Dubson~\cite{Dubson} and Kashiwara~\cite{Kindex}.
This states that given a constructible
complex of sheaves $\CF$, its Euler characteristic $\chi(X,\CF)$ is equal
to the intersection of Lagrangian cycles $CC(\CF)\cdot [df]$ where $df$ is the graph
of a sufficiently
generic function $f:X\to\R$.
More generally, given two constructible
complexes of sheaves $\CF_1,\CF_2$, a formula of MacPherson (see the introduction of~\cite{Gins}, 
the lectures notes~\cite{GMacP}, and a Floer-theoretic interpretation~\cite{KO1, KO2})
expresses the Euler characteristic of their tensor product in terms of the intersection
of their characteristic cycles
$$
\chi(X,\CF_1\stackrel{L}{\otimes} \CF_2)=CC(\CF_1)\cdot CC(\CF_2).
$$

The most direct influence on our main result
is the work of Ginsburg~\cite{Gins} (in the complex affine case) and Schmid-Vilonen~\cite{SV}
(in general)
on the functoriality of
the characteristic cycle construction. Thanks to their work, one knows how to calculate the 
characteristic cycle $CC(Ri_*\CF)$ of the direct image under an open embedding $i:U\hookrightarrow X$. (The functoriality of the characteristic cycle
under the other standard operations is explained by
Kashiwara-Schapira~\cite{KS}.) 
In the subanalytic context,
given an open subset $i:U\hookrightarrow X$, one can always choose a defining function
$m:X\to \R_{\geq 0}$
for the boundary $\del U\subset X$. By definition, $m$ is a nonnegative function whose zero
set is precisely $\del U\subset X$. 
With such a function in hand,
the formula for open embeddings is the limit of Lagrangian cycles
$$
CC(Ri_*\CF)=\lim_{\epsilon\to 0_+} (CC(\CF) + \epsilon \Gamma_{d\log m})
$$
where $\Gamma_{d\log m}\subset T^*U$ is the graph of the differential,
and the sum is set-theoretic.
The proof of our main result may be interpreted as a categorification of this formula.
We explain this in the next section.


\subsection{Summary}

To relate the constructible derived category $D_c(X)$ to the Fukaya category $Fuk(T^*X)$, we proceed in several steps, some topological and some categorical.

It is well-known that
usual notions of category theory are too restrictive a context for
dealing with the geometry of moduli spaces of pseudoholomorphic maps.
To be precise, $Fuk(T^* X)$ is not a usual category but rather an $A_\infty$-category.
Relations among compositions of morphisms are determined by the bubbling 
of pseudoholomorphic disks, and this is not associative but only homotopy associative.
The $A_\infty$-category formalism is a means to organize 
these homotopies (and the homotopies between the homotopies,
and so on). In particular, morphisms in an $A_\infty$-category are represented by
chain complexes to provide some homotopic flexibility.
When this is the only added wrinkle, so that compositions of such morphisms
are in fact associative, one arrives at the special case of a differential graded (dg) category.
To an $A_\infty$-category one can assign an ordinary (graded)
category by taking the cohomology of its morphism complexes.
This allows for the perspective that these notions only differ from that of an ordinary category
by providing more homotopic flexibility. (We collect some of
the basic notions of $A_\infty$-categories in Section~\ref{a-infinity background} below.)

The derived category $D_c(X)$ is the cohomology category
of a dg category $Sh(X)$ whose objects are constructible complexes of sheaves. 
The morphisms in $Sh(X)$ are defined by starting with the naive definition of morphisms
of complexes, and then passing to the dg quotient category
where quasi-isomorphisms are invertible.
Our first step in reaching $Fuk(T^*X)$ is to observe
that $Sh(X)$ is generated by its full subcategory consisting of 
standard objects associated to open submanifolds. 
In the subanalytic context,
given an open subset $U\subset X$, one can always choose a defining function
$m:X\to \R_{\geq 0}$
for the complement $X\setminus U$. By definition, $m$ is a nonnegative function whose zero
set is precisely $X\setminus U$. To keep track of the choice of such a function,
we define a dg category $Open(X)$ as follows. 
Its objects are pairs $(U,m)$ where $U\subset X$ is open,
and $m$ is as described.
Its morphisms are given by complexes of relative de Rham forms, and
are naturally quasi-isomorphic with those for the corresponding standard objects of $Sh(X)$.
In the language of dg categories, one can say that $Sh(X)$ is a triangulated envelope of $Open(X)$,
and that $D_c(X)$ is the derived category of both $Sh(X)$ and $Open(X)$.

The aim of our remaining constructions is to embed
 the $A_\infty$-category $Open(X)$ into the Fukaya $A_\infty$-category $Fuk(T^*X).$
It is simple to say where this $A_\infty$-functor takes objects of $Open(X)$.
To explain this, we introduce some notation in a slightly more general context.
Given a submanifold $Y\subset X$
and a defining function $m:X\to\R_{\geq 0}$
for the boundary $\del Y\subset X$, set
$f:X\setminus \partial Y \to \R$ to be the logarithm $f=\log m$, and define
the {standard Lagrangian} $L_{Y,f}\subset T^*X\vert_Y \subset T^*X$ to be
the fiberwise sum
$$
L_{Y,f}= T^*_Y X + \Gamma_{df}\vert_Y,
$$
where $\Gamma_{df}\subset T^*X\vert_{X\setminus \partial Y}$
denotes the graph of $df,$
and the sum is taken fiberwise in $T^*X\vert_Y.$
By construction,
$L_{Y,f}$ depends only on the restriction
of $m$ to $Y.$
In particular, 
for an open subset $U\subset X$, we could also take $m$ to be a defining function
for the complement $X\setminus U$. In this case, the definition simplifies so that
$
L_{U,f}$
is just the graph 
$
\Gamma_{df}
$
over $U$.

Now given an object $(U,m)$ of $Open(X)$, where $U\subset X$ is open,
and $m:X\to \R_{\geq 0}$ is a defining function for $X\setminus U$,
we send it to the standard Lagrangian $L_{U,f}\subset T^*X$, where $f:U\to\R$
is given by $f=\log m$. If $U$ is not all of $X$, this is a closed but noncompact Lagrangian
submanifold of $T^*X$. To properly obtain an object of $Fuk(T^*X)$,
we must endow $L_{U,f}$ with a brane structure. This consists of a 
grading (or
lifting of its squared phase) and relative pin structure. We check that standard Lagrangians
carry canonical brane structures with respect to canonical background classes.
We make $L_{U,f}$ an object
of $Fuk(T^*X)$ by equipping it with its canonical brane structure.

What is not immediately clear is what our $A_\infty$-functor should do with morphisms.
To answer this,
we first identify $Open(X)$ with an $A_\infty$-category $Mor(X)$
built out of the Morse theory of open subsets of $X$ equipped with
defining functions for their complements. The construction of $Mor(X)$
is a generalization of Fukaya's Morse $A_\infty$-category of a manifold.
As with $Open(X)$, the objects of $Mor(X)$
are pairs $(U,m)$, where $U\subset X$ is open,
and $m:X\to \R_{\geq 0}$ is a defining function for $X\setminus U$.
As usual, it is convenient to set $f=\log m$ as a function on $U$.
For a finite collection
of objects $(U_i,m_i)$ of $Mor(X)$ indexed by $i\in\Z/(d+1)\Z$, the morphisms and composition
maps among the objects
encode the moduli spaces of maps from trivalent trees into $X$ that take edges to
gradient lines of the functions $f_{j}-f_i$ with respect to some Riemannian
metric on $X$. For example, the morphism complexes
are generated by the critical points of  Morse functions on certain open subsets,
and the differentials are given by counting isolated gradient lines. 

There are several delicate aspects
to working out the details of this picture. As usual with such a construction, we must be sure that
the functions $f_i$ and the Riemannian metric are sufficiently generic to ensure
we have well-behaved moduli spaces. But in our situation, we must also be sure that the 
gradient vector fields of the differences $f_{i+1}-f_{i}$ are not wild at the boundaries of
their domains $U_i\cap U_{i+1}$. To accomplish this, we employ techniques
of stratification theory to move the boundaries and functions into a sufficiently
transverse arrangement. Then there will be an open, convex space of Riemannian metrics
such that the resulting moduli spaces are well-behaved. 
The upshot
is that we obtain an $A_\infty$-structure on $Mor(X)$ whose composition
maps count so-called gradient trees for Morse functions on certain open subsets of $X$.
Furthermore, an application of arguments of Kontsevich-Soibelman~\cite{kontsoib}
from homological perturbation theory
provides a quasi-equivalence 
$$
Open(X)\simeq Mor(X).
$$

Finally, we embed the Morse $A_\infty$-category $Mor(X)$ into the Fukaya $A_\infty$-category
$Fuk(T^*X)$ as follows. 
Let $(U_i,m_i)$ be a collection of objects of $Mor(X)$ indexed by $i\in\Z/(d+1)\Z$, 
and let $L_{U_i,f_i}$ be the corresponding collection of standard branes of $Fuk(T^*X)$
where as usual $f_i=\log m_i$.
After carefully perturbing the objects,
we check that the moduli spaces of gradient trees for the former collection may be identified
with the moduli spaces of pseudoholomorphic polygons for 
the latter.
When all of the open sets $U_i$ are the entire manifold $X$,
this is a theorem of Fukaya-Oh~\cite{FO}.
These authors have identified the Morse $A_\infty$-category 
of the manifold $X$ and the Fukaya $A_\infty$-category
of graphs in $T^*X$.
To generalize this to arbitrary open sets, we employ the following strategy.
First, using area bounds,
we check that all pseudoholomophic maps with boundary on our standard branes
in fact have boundary in a prescribed
region. Next, we dilate our standard branes
so that the theorem of Fukaya-Oh identifies 
the relevant moduli subspaces.
Finally, we check that the homogeneity of the area bounds
guarantees that no critical event occurs during the dilation.
Thus we obtain an $A_\infty$-embedding
$$
Mor(X)\hookrightarrow Fuk(T^*X).
$$

Putting together the above functors gives a quasi-embedding of the $A_\infty$-category $Sh(X)$
of constructible complexes of sheaves on $X$
into the $A_\infty$-category 
of twisted complexes $\tw Fuk(T^*X)$ in the Fukaya category
of $T^*X$. Taking the underlying
cohomology categories gives an embedding of the corresponding derived categories.

For future applications, it is useful to know where the embedding takes other objects 
and morphisms. In particular, we would like to know not only where it takes standard sheaves
on open submanifolds, but also standard sheaves on arbitrary submanifolds. One approach
to this problem is to express standard sheaves on arbitrary submanifolds in terms of standard sheaves
on open submanifolds, and then to check what the relevant distinguished
triangles of constructible sheaves look like under the embedding. This requires
identifying certain cones in the Fukaya category with symplectic surgeries. 
Rather than taking this 
route,
we will instead show in the final section that we may explicitly extend the domain of the embedding to include
standard sheaves on arbitrary submanifolds and morphisms between them. 
 This has the added
value that given constructible sheaves on a stratification, it obviates the need to further
refine the stratification in order to construct the embedding: one may use the standard
sheaves themselves as a generating set.

Consider the standard sheaf $Ri_*\CL_Y$
associated to a local system $\CL_Y$ on an arbitrary submanifold $i:Y\hookrightarrow X$.
Suppose that we are given a defining function $m:X\to \R_{\geq 0}$ for the boundary $\del Y\subset X$.
Recall that we define the standard Lagrangian $L_{Y,m}\subset T^*X$ to be the fiberwise sum
$$
L_{Y,f} =T^*_Y X + \Gamma_{df}
$$
where $T^*X_Y\subset T^*X$ is the conormal bundle to $Y$, and $\Gamma_{df}\subset T^*X$
is the graph of the differential of $f=\log m$.
We write $L_{Y,f,\CL_Y}$ for the corresponding standard object of $Fuk(T^*X)$
given by $L_{Y,f}$ equipped with its canonical brane structure
and the pullback of the flat vector bundle $\CL_Y\otimes or_X\otimes or_Y^{-1}$,
where $or_X, or_Y$ denote the orientation bundles of $X,Y$ respectively.
The main consequence of the final section is the following.

\begin{thm}
Under the embedding $D_c(X)\hookrightarrow D Fuk(T^*X)$, the image of
the standard sheaf $Ri_*\CL_Y$ is canonically isomorphic to the standard brane $L_{Y,f,\CL_Y}$.
\end{thm}


\subsection{Mirror symmetry}

The connection of this current work to mirror symmetry is
somewhat speculative, though
several appearances of constructible sheaves
in the context of mirror symmetry deserve mention.

First, the announced results of Bondal and Bondal-Ruan~\cite{BR}
relate the
derived categories of coherent sheaves on toric
Fano varieties with the Fukaya-Seidel category
on the Landau-Ginzburg side. Their method is to establish equivalences of both with
the derived category of constructible sheaves on a torus with respect
to a specific (non-Whitney) stratification determined by the superpotential.
One can view the result of Bondal-Ruan
from the perspective developed in this paper by
identifying $(\mathbb C^*)^n$
with $T^*((S^1)^n).$

Second, Kapustin-Witten~\cite{KW} place the geometric Langlands program
in the context of topological quantum field theory. In particular,
they relate the harmonic analysis of the
geometric Langlands program to mirror symmetry by equating
Hecke operators on $\CD$-modules with 't Hooft operators
acting on branes.
In this setting, one may interpret the results of this paper as lending some mathematical evidence
to this physical perspective. For example, according to Kapustin-Witten~\cite{KW}, 
given a generic eigen-brane for the 't~Hooft operators,
there is a corresponding regular, holonomic Hecke eigen-$\CD$-module.
One might hope to provide an explicit construction of the eigen-$\CD$-module by
first identifying the eigen-brane as the microlocalization of some constructible sheaf,
and then applying the Riemann-Hilbert correspondence.

Third, braid group actions have been an active area of interest especially
in the context of branes in the cotangent bundle of flag varieties $\CB$.
In the case of coherent sheaves, braid group actions on $D^b_{coh}(T^*\CB)$ 
have been studied by many authors (see for example Seidel-Thomas~\cite{SeidelThomas}). 
One may use the results of this paper
to construct the corresponding actions in the symplectic context.
Namely, under the embedding of this paper,
the kernels giving the usual braid group action
on the constructible derived category $D_c(\CB)$ (see for example Rouquier~\cite{RouquierBraid})
induce a corresponding action on $DFuk(T^*\CB)$.

Fourth, the work of Kontsevich-Soibelman
\cite{kontsoib2}
and Gross-Siebert \cite{GS} paints the
large complex structure limit of a Calabi-Yau $n$-fold as a
collapse into a real $n$-fold with integral affine structure
and a Monge-Amp\`ere metric.  The complex $n$-fold
is recovered from the limit manifold as a quotient of the
tangent (or cotangent) bundle by the lattice of integer tangent
vectors.  It is intriguing to imagine
a quotient construction creating a torus fibration from
the cotangent bundle.

Finally, it would be interesting to understand whether our result is the local picture
of a relationship that holds more generally for compact symplectic manifolds.
One may consider modules over the deformation quantization
as a global analogue of constructible sheaves. (See for example Kontsevich~\cite{kontquant}, Kashiwara~\cite{Kquant}, or Polesello-Schapira ~\cite{PS}.)
Clear comparisons can then be made between such modules and the Fukaya category.
There is great interest in understanding more precisely how to interpolate between
the local nature of the modules and the global nature of the Fukaya category.

\begin{ack}
We have benefited greatly from discussions with Melissa Liu, Paul Seidel, 
and Chris Woodward. We would also like to thank Yong-Geun Oh for comments
about the context of our work and its exposition. Finally, we are grateful to
an anonymous referee whose comments have led to significant improvements in the paper.

The work of D.\ N. is supported in part by NSF grant DMS--0600909 and DARPA. 
The work of E.\ Z. is supported in part by NSF grant
DMS--0405859.
Any opinions, findings and conclusions or recommendations expressed
in this material are those of the authors and do not necessarily reflect the
views of the National Science Foundation (NSF).
\end{ack}


\section{$A_\infty$-categories}\label{a-infinity background}
We collect here standard material concerning $A_\infty$-categories,
dg categories, and triangulated categories. Our reference is Chapter 1 of
Seidel's book \cite{Seidel}.


\subsection{Preliminaries}
 Our aim here is not to recall complete
definitions, but only to establish notation. 

Let $\CA$ be a (not necessarily unital) $A_\infty$-category with set of objects $\ob\CA$, $\Z$-graded
vector space of morphisms $\hom_\CA(X_0,X_1)$,
and composition maps
$$
\mu^d_\CA: \hom_\CA(X_0,X_1)\otimes \cdots \otimes \hom_\CA(X_{d-1},X_d)
\to \hom_\CA(X_0,X_d)[2-d],
\mbox{ for $d\geq 1$}.
$$
A dg category is an $A_\infty$-category $\CA$ whose higher composition maps
$\mu_\CA^d$, for $d\geq 3$ are equal to zero.

Let $H(\CA)$ denote the $\Z$-graded
cohomological category of $\CA$ with set of objects $\ob H(\CA)=\ob\CA$, and $\Z$-graded
vector space of morphisms 
$$
\hom_{H(\CA)}(X_0,X_1) = H(\hom_\CA(X_0,X_1),\mu^1_\CA).
$$
Let $H^0(\CA)$ denote the ungraded cohomological category with 
set of objects $\ob H^0(\CA)=\ob\CA$, and 
vector space of morphisms 
$$
\hom_{H^0(\CA)}(X_0,X_1) = H^0(\hom_\CA(X_0,X_1),\mu^1_\CA).
$$
An $A_\infty$-category is said to be {\em cohomologically unital} or {\em c-unital}
if $H(\CA)$ is unital.

Let $\CF:\CA\to\CB$ be an $A_\infty$-functor between $A_\infty$-categories
with map on objects $\CF:\ob\CA\to\ob\CB$, and morphism maps
$$
\CF^d: \hom_\CA(X_0,X_1)\otimes \cdots \otimes \hom_\CA(X_{d-1},X_d)
\to \hom_\CB(\CF X_0,\CF X_d)[1-d],
\mbox{ for $d\geq 1$}.
$$
An $A_\infty$-functor is said to be {\em c-unital} if $H(\CF)$ is unital.

Throughout what follows, we assume that all $A_\infty$-categories are c-unital,
and all $A_\infty$-functors are c-unital.
We say that an $A_\infty$-functor $\CF$ is 
a {\em quasi-equivalence} if the induced functor $H(\CF)$ is an equivalence.
We say that $\CF$ is a {\em quasi-embedding} if $H(\CF)$ is full and faithful.


\subsection{$A_\infty$-modules}

Let $\ch$ denote the dg category of chain complexes, considered as an $A_\infty$-category.

Given an $A_\infty$-category $\CA$, 
an $A_\infty$-module over $\CA$ is an $A_\infty$-functor 
$\CA^{opp}\to \ch.$
Let $\modules(\CA)$ 
denote the  $A_\infty$-category of $A_\infty$-modules over $\CA$.

The functor category 
$\modules(\CA)$ inherits much of the structure of
the target category $\ch$. For example, $\modules(\CA)$ is a dg category,
and its cohomological category $H^0(\modules(\CA))$ is a triangulated category.
In particular,
we have the obvious shift functor $S$ on modules and the cohomological notion of exact triangle
of modules. Note that the shift functor may be recovered by taking the cone of the zero morphism
to the trivial module, or to the cone of the identity morphism of any module.

For any object $Y\in\ob\CA$, we have the $A_\infty$-module $\CY(X)=\hom_\CA(X,Y)$
with $\mu^d_\CY=\mu^d_\CA$.
This provides an $A_\infty$-Yoneda embedding
$
\CJ:\CA\to \modules(\CA)
$
which is cohomologically full and faithful.
Since the ambient category $\modules(\CA)$ is a dg category,
the image $\CJ(\CA)$ of the Yoneda embedding is as well.
Thus each $A_\infty$-category $\CA$
is canonically quasi-equivalent to a dg category $\CJ(\CA)$.


\subsection{Triangulated $A_\infty$-categories}

Given an $A_\infty$-category $\CA$, 
an {\em exact triangle} in $H(\CA)$ is defined to be any diagram in $H(A)$
which becomes isomorphic to an exact triangle of $H(\modules(\CA))$ under the Yoneda embedding.
A {\em shift} $SX$ of an object $X$ is any object which becomes isomorphic to the shift 
in $H(\modules(\CA))$ under the Yoneda embedding.

A nonempty $A_\infty$-category $\CA$ is said to be {\em triangulated}
if the following hold:

\begin{enumerate}
\item Every morphism in $H^0(\CA)$ can be completed to an exact triangle in $H(\CA)$.
In particular, every object $X$ has a shift $SX$.

\item For each object $X$, there is an object $\tilde X$
such that $S\tilde X\simeq X$ in $H^0(\CA)$.

\end{enumerate}

If $\CA$ is a triangulated $A_\infty$-category, then $H^0(\CA)$ is a triangulated
category in the usual sense. Furthermore, if $\CF:\CA\to\CB$ is an $A_\infty$-functor
between triangulated $A_\infty$-categories, then $H^0(\CF)$ is an exact functor.

Let $\CA$ be a full $A_\infty$-subcategory of a triangulated $A_\infty$-category $\CB$.
The triangulated $A_\infty$-subcategory of $\CB$ {\em generated} by $\CA$ is the smallest
full subcategory $\ol \CA$ that contains $\CA$, is closed under cohomological isomorphism,
and is itself triangulated.

A {\em triangulated envelope} of a nonempty $A_\infty$-category $\CA$
is a pair $(\ol\CA,\CF)$ consisting of a triangulated $A_\infty$-category $\ol\CA$,
and a cohomologically full and faithful functor $\CF:\CA\to\ol\CA$ such that
the objects in the image of $\CF$ generate $\ol\CA$.
The triangulated category $H^0(\ol\CA)$ is independent of the choice of envelope
up to exact equivalence. 
It is sometimes called the {\em derived category} of $\CA$ and denoted by $D(\CA)$,
but we will sometimes reserve this to mean a localized version of $H^0(\ol\CA)$.
Thus when we use the term derived category and the notation $D(\CA)$,
we will be explicit about what is intended.


\subsection{Twisted complexes}
There are two standard constructions of a triangulated envelope: (i)
the full subcategory of $\modules(\CA)$ generated by the image of the Yoneda embedding, 
and (ii) the $A_\infty$-category of twisted complexes $\tw\CA$. 

In this paper, we adopt the approach of twisted complexes.
The explicit construction of $\tw\CA$ will play no role, only the following formal properties.

First, $\tw\CA$ is a triangulated $A_\infty$-category.
There is a canonical 
$A_\infty$-functor 
$\iota:\CA\to\tw\CA$ such that $\iota$ is injective on objects, on morphisms we have
$$\hom_\CA(X_0,X_1)=\hom_{\tw\CA}(\iota X_0,\iota X_1),$$
and the composition maps $\mu_\CA^d$ and $\mu_{\tw\CA}^d$ coincide for objects of $\CA$
and their images under $\iota$.
In short, we may identify $\CA$ with its image under $\iota$.

We also have the following.

\begin{enumerate}
\item If $\CA$ is c-unital, then $\tw\CA$ and $\iota$ are as well.
\item $\tw\CA$ is generated by $\CA$.
\end{enumerate}

Furthermore, any $A_\infty$-functor $\CF:\CA\to\CB$ extends to an $A_\infty$-functor 
$\tw\CF:\tw\CA\to\tw\CB$ satisfying the following.

\begin{enumerate}
\item If $\CF$ is c-unital, then $\tw\CF$ is as well.
\item If $\CF$ is cohomologically full and faithful, then $\tw\CF$ is as well.
\item If $\CF$ is a quasi-equivalence, then $\tw\CF$ is as well.

\end{enumerate}


\subsection{Homological perturbation theory}
\label{spt}
We recall here the general picture of homological perturbation theory as
summarized by Seidel \cite{Seidel}. 

Let $\CB$ be an $A_\infty$-category. Suppose that for each pair of objects $(X_0,X_1)$,
we have a chain complex $(\hom_\CA(X_0,X_1),\mu^1_\CA)$, chain maps
$$
\CF^1:\hom_\CA(X_0,X_1)\to\hom_\CB(X_0,X_1),\qquad
\CG^1:\hom_\CB(X_0,X_1)\to\hom_\CA(X_0,X_1)
$$
of degree $0$, and a linear map
$$
T^1: \hom_\CB(X_0,X_1)\to \hom_\CB(X_0,X_1)
$$
of degree $-1$ such that
$$
\CF^1\circ \CG^1-id=
\mu^1_B T^1+T^1 \mu_\CB^1. 
$$

In the preceding set-up, the subscript $\CA$ is simply suggestive notation.
The main statement of homological perturbation theory is that 
there is an explicit construction of an
$A_\infty$-category $\CA$ with objects $\ob\CA=\ob\CB$, and
morphism complexes the given $(\hom_\CA(X_0,X_1),\mu^1_\CA)$. 
Furthermore, there are $A_\infty$-functors $\CF:\CA\to\CB,$ $\CG:\CB\to\CA$ which are the identity
on objects, and have first-order terms the given $\CF^1,\CG^1$. Finally,
there is a homotopy between $\CF\circ\CG$ and $id_\CB$ which starts with the given $T^1$.

We will use the special case of this construction when $\CG^1$ is an idempotent $\pi^1$,
and $\CF^1$ is the inclusion $i^1$ of the image of $\pi^1$.
In other words, for each pair of objects $(X_0,X_1)$, we have a chain map
$$
\pi^1:\hom_\CB(X_0,X_1)\to\hom_{\CB}(X_0,X_1)
$$
of degree $0$ such that $\pi^1\circ\pi^1=\pi^1$,
and a linear map
$$
T^1: \hom_\CB(X_0,X_1)\to \hom_\CB(X_0,X_1)
$$
of degree $-1$ such that
$$
i^1\circ \pi^1-id=
\mu^1_B T^1+T^1 \mu_\CB^1. 
$$
In this case, if we take 
$$
\hom_\CA(X_0,X_1)= \pi^1(\hom_\CB(X_0,X_1)),
$$
then the resulting $A_\infty$-functors $i:\CA\to\CB$, $\pi:\CB\to\CA$
are quasi-equivalences.


\section{Analytic-geometric categories}
\label{subanalytic}

When working with sheaves on a manifold $X$, it is often useful
if not indispensable to restrict to subsets of $X$ that have strong finiteness properties.
In this section, we collect basic material from the theory of subanalytic sets 
that
plays a role in what follows. All of the results and arguments that we use hold
in the context of analytic-geometric categories.
Since this seems to be a natural level of generality, we adopt it as our framework.
What follows is a brief summary of relevant results from van den Dries-Miller~\cite{vdDM}.
For a discussion of subanalytic sets alone, see Bierstone-Milman~\cite{BM}.  The
reader may prefer to ignore the generality of analytic-geometric categories 
and consider all discussion to take place in the subanalytic category.

Throughout what follows, all manifolds are assumed to be real analytic unless otherwise specified.
 
\subsection{Basic definitions}
To give an {\em analytic-geometric category } $\CC$ is to equip
each manifold $M$ with a collection $\CC(M)$ of subsets of $M$
satisfying the following properties:

\begin{enumerate}
\item $\CC(M)$ is a Boolean algebra of subsets with $M\in\CC(M)$.

\item If $A\in\CC(M)$, then $A\times\R\in\CC(M\times \R)$.

\item If $f:M\to N$ is a proper analytic map and $A\in\CC(M)$, then $f(A)\in\CC(N)$.

\item If $A\subset M$, and $(U_i)_{i\in I}$ is an open covering of $M$, then $A\in\CC(M)$
if and only if $A\cap U_i\in\CC(U_i)$, for all $i\in I$.

\item For every bounded set $A\in \CC(\R)$, the boundary $\del A$  is finite. 
\end{enumerate}

Given the above data, one defines a category $\CC$ as follows.
An object is a pair $(A,M)$ with $M$ a manifold, and $A\in\CC(M)$.
A morphism $(A,M)\to (B,N)$ is a continuous map $f:A\to B$ whose graph
$\Gamma(F)\subset A\times B$ belongs to $\CC(M\times N)$.
Objects of $\CC$ are called $\CC$-sets, and morphisms are called $\CC$-maps.
When the codomain of a map is $\R$, we refer to it as a function.

The basic example of an analytic-geometric category is the subanalytic category $\CC_{an}$
of subanalytic sets and continuous maps with subanalytic graphs.
For any analytic-geometric category $\CC$, the subanalytic subsets of
any manifold $M$ belong to $\CC(M)$.


\subsection{Background results}
Most of the fundamental results about subanalytic sets hold in any analytic-geometric category 
(although it is unknown whether the uniformization and rectilinearization properties
of subanalytic sets have analogues).
We limit our discussion here to include only those results which we use.


\subsubsection{Derivatives}

Given a manifold $M$, the tangent bundle $TM$ and cotangent bundle $T^*M$ are also manifolds.
Given a $C^1$ submanifold $A\subset M$, let $TA\subset TM$ denote its tangent bundle, and
$T^*_A M\subset T^*M$ its conormal bundle.

\begin{lem}\label{derivs}
If $A\in \CC(M)$ is a $C^1$ submanifold of $M$, then $TA\in\CC(TM)$ and $T^*_A M\in\CC(T^*M)$.
If $f:M\to N$ is a $\CC$-map of class $C^1$, its differential $Tf:TM\to TN$ is a $\CC$-map.
\end{lem}
\subsubsection{Whitney stratifications}

Let $X,Y$ be $C^1$ submanifolds of a manifold $M$, and let $x\in X$. 
The triple $(Y, X,x)$ is said to satisfy {\em Whitney's condition} 
if 
given any sequences of points $x_i\in X$ and 
$y_i\in Y$ each converging to $x$, such that in some local coordinate chart
the secant lines $\ell_i=\ol{x_i y_i}$
converge to some line $\ell$ and the tangent planes $T_{y_i} Y$
converge to some plane $\tau$, we have
$
\ell\subset \tau.
$
The pair $(Y,X)$ is said to satisfy {Whitney's condition} if  for all $x\in X$,
the
triples $(Y,X,x)$ satisfy
the condition.

A {\em $C^p$ stratification} of a manifold $M$
consists of a locally finite collection $\CS=\{S_\alpha\}$ of locally closed $C^p$ submanifolds
$S_\alpha\subset M$ called {\em strata} satisfying

\begin{enumerate}
\item
(covering) $X=\bigcup_\alpha S_\alpha$, 
\item
(pairwise disjoint) $S_\alpha \cap S_\beta =\emptyset,$ for $\alpha\not =\beta$, 
\item
(axiom of frontier) 
$\ol {S_\alpha} \cap S_\beta \not =\emptyset$ if and only if $S_\beta \subset \ol{S_\alpha}$.
\end{enumerate}

A $C^p$ stratification $\CS=\{S_\alpha\}$ of $M$ is called a {\em $C^p$ Whitney stratification} if 
$(S_\alpha,S_\beta)$ satisfies Whitney's condition for all $\alpha,\beta$.

A stratification $\CS$ of $M$ is said to be {\em compatible} 
with a collection $\CA$ of subsets of $M$ if 
$S\cap A\not =\emptyset$ implies $S\subset A$,
for all $S\in\CS, A\in\CA$.

Given a map $f:M\to N$ between manifolds, a {\em $C^p$ Whitney stratification}
of $f$ is a pair $(\CS,\CT)$ where $\CS$ and $\CT$ are $C^p$ Whitney stratifications
of $M$ and $N$ respectively such that for each $S\in \CS$, the map $f|_S:S\to N$
is a $C^p$ submersion with $f(S)\in\CT$.

\begin{prop}\label{whitney}
Let $X\in\CC(M)$ be closed, and $p$ a positive integer.
\begin{enumerate}
\item For every locally finite collection $\CA\subset \CC(M)$, there is a $C^p$ Whitney
stratification $\CS\subset \CC(M)$ of $M$ which is compatible with $\CA$
and has connected strata.

\item Let $f:M\to N$ be a proper $\CC$-map, and let $\CA\subset\CC(M),\CB\subset \CC(N)$
be locally finite collections. Then there is a $C^p$ Whitney stratification $(\CS,\CT)$ of $f$
with connected strata such that $\CS\subset \CC(M)$ is compatible with $\CA$ and $\CT\subset \CC(N)$
is compatible with $\CB$.
\end{enumerate}
\end{prop}

\begin{rmk}\label{cell}
In Proposition~\ref{whitney}(1),
we may find a stratification $\CS\subset \CC(M)$ of $M$ such that each stratum $S\in \CS$
is a $(C^p,\CC)$ cell in $M$.
In particular, each $S\in \CS$ will be $C^p$ diffeomorphic to $\R^d$, where $d=\dim S$. 
\end{rmk}

To a Whitney stratification $\CS=\{S_\alpha\}$ of $M$, we associate the conical
set $\Lambda_\CS \subset T^*M$ given by the union
$$
\Lambda_\CS =\cup_\alpha T^*_{S_\alpha} M.
$$
By Lemma~\ref{derivs}, if $\CS\subset \CC(M)$, then $\Lambda_\CS\in \CC(T^*X)$.
If a stratification $\CS_1$ of $M$ is compatible with another $\CS_2$, then 
$\Lambda_{\CS_2}\subset \Lambda_{\CS_1}$.

For a function $f:X\to \R$, 
we say that $x\in X$ is a {\em $\Lambda_\CS$-critical point} of $f$ if we have $df(x)\in\Lambda_\CS$.
We say that $r\in \R$ is a {\em $\Lambda_\CS$-critical value} of $f$ 
if there is
a $\Lambda_\CS$--critical point $x\in X$ such that $r=f(x)$, 
otherwise we say that $r$ is a {\em $\Lambda_\CS$-regular value}.

\begin{lem}\label{sard}
If $f$ is a $\CC$-map,
then the $\Lambda_\CS$-critical values of $f$ form a discrete subset of $\R$.
\end{lem}

\begin{proof}
The $\Lambda_\CS$-regular values are a dense open $\CC$-subset of $\R$.
\end{proof}


\subsubsection{Curve selection lemma}

\begin{prop}\label{csl}
If $A\in\CC(M)$, and $x\in \ol A\setminus A$, then there is a $\CC$-map $\gamma:[0,1)\to M$
such that $\gamma(0,1)\subset A$, and $\gamma(0)=x$. Furthermore, if $p$ is a positive integer,
$\gamma$ can be chosen to be injective and of class $C^p$.
\end{prop}


\subsubsection{Defining functions}\label{defining}
For a function $m:M\to\R$, let $Z(m)=\{x\in M|m(x)=0\}$. 
Given a subset $A\subset M$, we call any function  $m:M\to\R_{\geq 0}$ with $Z(m)=A$
a {\em defining function} for $A$.
The space of defining functions for $A$ is convex.

\begin{prop}
If $A\in\CC(M)$ is closed and $p$ is a positive integer, then there is a $\CC$-map
$m:M\to\R$ of class $C^p$ with $A=Z(m)$.
\end{prop}

Using defining functions, one can construct bump functions as follows. For
$A_1,A_2\in\CC(M)$ disjoint and closed with defining functions $m_1, m_2$ respectively,
define $g_i=m_i^2/(1+m_i^2)$, for $i=1,2$. Then the function $b=(g_1 +g_1 g_2)/(g_1+g_2)$
satisfies $Z(b)=A_1$, $Z(b-1)=A_2$, and $0\leq b\leq 1$.
The existence of bump functions implies the existence of partitions of unity
and constructions which devolve from them.

\begin{rmk} 
In simple situations, one can explicitly produce $C^p$ bump functions with compact support
 such as
$$
b(x) = \begin{cases} (1 - x^2)^{p+1}, \quad x^2<1,\\ \quad\; 0,
\qquad \mbox{otherwise}.\end{cases}
$$
\end{rmk}


\section{Constructible sheaves}\label{constr}

Let $X$ be a real analytic manifold. All subsets of $X$ are assumed to belong
to some fixed analytic-geometric category unless otherwise specified.

Let $\C_X$ be the sheaf of locally-constant complex-valued functions on $X$.
By a sheaf on $X$, we will always mean a sheaf of $\C_X$-modules.
A sheaf $\CF$ is said to be {\em constructible} if there exists a Whitney
stratification of $X$ such that the restriction of $\CF$ to each stratum is locally-constant
and finitely-generated.

We define
the localized triangulated dg category $Sh(X)$ of complexes of sheaves with bounded constructible cohomology
as follows.
First, we have the naive triangulated  dg category $Sh_{naive}(X)$ whose objects are
complexes of sheaves with bounded constructible cohomology, and whose morphisms are
the usual complexes of morphisms. Then, we take $Sh(X)$ to be the dg quotient
of $Sh_{naive}(X)$ with respect to the subcategory $\CN$ of acyclic objects~\cite{Keller}.
As explained in~\cite{Drinfeld}, this can be achieved by simply 
adding a homotopy between zero and the identity to the endomorphism complex of each object of $\CN$.
The ungraded cohomological category $H^0 (Sh(X))$ is the usual bounded constructible derived category $D_c(X).$

We have 
the six standard derived functors of Grothendieck
$f_*$, $f^*$, $f_!$, $f^!$, $\otimes$ and $\intHom$. 
We similarly have the Verdier duality functor $\CD$. 
Note that we only consider derived functors though the notation does not make this explicit.
We also refer to objects of $Sh(X)$ as sheaves though they
are properly complexes of sheaves.


\subsection{Standard objects}

The most accessible objects of $Sh(X)$ are the so-called standard and costandard sheaves
of submanifolds.
To be precise, let $i:Y\hookrightarrow X$ be the inclusion of a submanifold 
(with its subspace topology) with closure $\ol Y\subset X$
and boundary $\del Y= \ol Y\setminus Y\subset X$. Note that the boundary could be a singular
subset.
For a local system $\CL_Y$ on $Y$,
we call the sheaf $i_*\CL_Y$ a {\em standard
object},
and the sheaf $i_!\CL_Y$ a {\em costandard object}. 
The terminology reflects that Verdier duality intertwines the two extensions
$$
\CD_X(i_!\CL_Y)\simeq i_*\CD_Y(\CL_Y).
$$
Suppose $U\subset X$ is an open set. The complex of sections of $i_*\CL_Y$ over $U$ 
is quasi-isomorphic to the complex of $\CL_Y$-valued singular cochains
$$
\Gamma(U,i_*\CL_Y)\simeq C^*(U\cap Y,\CL_Y).
$$
Similarly,
the complex of sections of $i_!\C_Y$ over $U$ 
is quasi-isomorphic to the complex of $\CL_Y$-valued relative singular cochains
$$
\Gamma(U,i_!\C_Y)\simeq C^*(U\cap \ol Y,U\cap \del Y,\CL_Y).
$$

\subsection{Standard triangles}
Fix a subset $Z\subset X$.
Let $i:U\hookrightarrow Z$ be the inclusion of an open set, and 
$j:Y=Z\setminus U\hookrightarrow Z$ the inclusion of its complement.
Note that since $U$ is open and $Y$ is closed in $Z,$
we have identities $i^!\simeq i^*$ and $j_!\simeq j_*$.
For any sheaf $\CF$ on $Z,$ we have the standard distinguished triangles
$$
j_* j^! \CF\to \CF\to i_* i^*\CF\stackrel{[1]}{\to}
\qquad
i_* i^! \CF\to \CF\to j_* j^*\CF\stackrel{[1]}{\to}
$$
For example, if we take $\CF$ to be $\C_Z$ and take the cohomology of global sections, we obtain 
from these two the standard
long exact sequences
$$
H^*(Z,U)\to H^*(Z)\to H^*(U)\stackrel{[1]}{\to}
\qquad
H^*(Z,Y)\to H^*(Z)\to H^*(Y)\stackrel{[1]}{\to}
$$

We also have distinguished triangles associated to truncation functors.
Let $\tau_{\leq \ell}$ be the functor which assigns to a complex $\CF$
the truncated complex 
$$
\cdots \to \CF^{\ell-1} \to \ker(d_\ell) \to 0 \to \cdots
$$
The natural map $\tau_{\leq \ell}\CF\to\CF$ induces an isomorphism on cohomology sheaves 
in degrees less than or equal to $\ell.$
Let $\tau_{>\ell}$ be the functor which assigns to a complex $\CF$
the truncated complex 
$$
\cdots \to 0  \to \on{im}(d_\ell) \to \CF^{\ell+1}\to \cdots
$$
The natural map $\CF\to\tau_{>\ell}\CF$ induces an isomorphism on cohomology sheaves 
in degrees greater than $\ell.$
We have a distinguished triangle
$$
\tau_{\leq \ell}\CF\to\CF\to\tau_{>\ell}\CF\stackrel{[1]}{\to}
$$

\subsection{Standard objects generate}

\begin{prop}
Any object of $Sh(X)$ is isomorphic to one obtained from shifts of standard objects
by iteratively forming cones.
The same is true for costandard objects. 
\end{prop}

\begin{proof}
Let $\CF$ be an object of $Sh(X)$. 
Fix a stratification $\CS$ 
of $X$ such that the cohomology sheaves of $\CF$ are constructible with respect to $\CS$.

We prove the first assertion (the second is similar, or follows by Verdier duality).
The proof is an induction on the strata, beginning with the open strata.
Let $i_k:\CS_k\to X$ be the inclusion of the union of the strata of dimension equal to $k$,
and let $j_{< k}:\CS_{< k}\to X$ be the inclusion of the union of the strata of dimension less than $k$.

Suppose $X$ has dimension equal to $n$. Then for the sheaf $\CF$,
we have a distinguished triangle
$$
j_{< n*} j_{< n}^! \CF\to \CF\to i_{n*} i_n^*\CF\stackrel{[1]}{\to}
$$
Using truncation functors, we may express the sheaf $\CF_n= i_{n*} i_n^*\CF$
by iteratively forming cones of shifted standard objects associated to the strata $\CS_n$.
By construction, the sheaf $\CF_{<n}=j_{< n*} j_{< n}^! \CF$ is supported on $\CS_{<n}$.

Next we have the distinguished triangle
$$
j_{< n-1*} j_{< n-1}^! \CF_{<n}\to \CF_{<n}\to i_{n-1*} i_{n-1}^*\CF_{<n}\stackrel{[1]}{\to}
$$
Again, using truncation functors, we may express the sheaf
$\CF_{n-1}=i_{n-1*} i_{n-1}^*\CF_{<n}$
by iteratively forming cones of shifted standard objects associated to the strata $\CS_{n-1}$.
By construction,
the sheaf $\CF_{<n-1}=j_{< n-1*} j_{< n-1}^! \CF_{<n}$ is supported on $\CS_{<n-1}$.

And so on. In the end, we see that $\CF$ may be expressed by iteratively forming
cones of shifted standard objects.
\end{proof}

We have the following strengthening of the proposition. 

\begin{prop} \label{generators}
Any object of $Sh(X)$ is isomorphic to one obtained from shifts of constant standard objects
$i_*\C_U$ for open submanifolds $i:U\hookrightarrow X$
by iteratively forming cones.
The same is true for constant costandard objects $i_!\C_U$. 
\end{prop}

\begin{proof} Let $\CF$ be an object of $Sh(X)$. 
Choose a stratification $\CT$
of $X$ such that the cohomology sheaves of $\CF$ are constructible with respect to $\CT$,
and the strata of $\CT$ are cells (see Remark~\ref{cell}).
By the previous proposition, and since the strata are cells,
we need to show that for a stratum $j:T\to X$,
we can realize the standard object $j_*\C_T$. 

Let $Star(T)$ be the union of all the strata of $\CT$ which
contain $T$ in their closures, and let $s:Star(T)\to X$ denote its inclusion. 
Let $Star'(T)$ be the complement of $T$ in $Star(T)$, and let 
$s':Star'(T)\to X$ denote its inclusion. Both $Star(T)$ and $Star'(T)$
are open submanifolds. We have the distinguished triangle
$$
j_*j^!\C_{Star(T)}\to s_* \C_{Star(T)}\to s'_*\C_{Star'(T)}\stackrel{[1]}{\to}
$$
Since $T\hookrightarrow Star(T)$ is the inclusion of an orientable submanifold, 
$j^!\C_{Star(T)}$ is isomorphic to a shift of $\C_{T}$, and the assertion is proved.
\end{proof}


\subsection{Open submanifolds}\label{open}
From here on, we focus on standard objects rather than costandard objects, though
there is no reason to prefer one over the other. Furthermore, we work with standard
objects for {open} submanifolds. By Proposition~\ref{generators}, such objects
generate the entire category $Sh(X)$. Thus it suffices to work with them
in proving our main theorem. While this will simplify many constructions,
there is a price to pay. First, we will lose concrete touch with other objects and only
understand what is happening with them in an abstract sense -- this has implications
for applications of our main result. Second, there
are contexts in which the arguments of Proposition~\ref{generators} are
not really acceptable -- often we are presented with a fixed stratification and
would prefer not to subdivide it further. To remedy both of these points, 
we have included a discussion in Section~\ref{arbitrary} explaining how to generalize
our arguments to deal with all standard objects, not only those 
for open submanifolds. 
In what follows, we also consider standard
objects with trivial coefficients and leave the case of arbitrary local systems to the reader.
This is a purely expositional choice, and the reader will have no trouble extending our
arguments. In any case, technically speaking, Proposition~\ref{generators} also obviates the need to consider arbitrary
local systems.

For an open subset $U\subset X$, let $\Omega^k(U)$ denote the space of differential $k$-forms
on $U$, and let $(\Omega(U),d)$ denote the deRham complex. 
Define the {\em support} of a $k$-form $\omega\in\Omega^k(U)$
to be the smallest closed subset $\supp(\omega)\subset X$ such that 
$$\omega|_{U\cap (X\setminus \supp(\omega))}=0.
$$
Now consider a pair $(U,V)$ where $V\subset U\subset X$ such that $U\setminus V$ is open.
Let $(\Omega(U,V),d)$
denote the relative deRham complex of differential forms on $U\setminus V$ whose support lies 
in $X\setminus V$. 
The complex 
$(\Omega(U,V),d)$ calculates the relative cohomology $H^*(U, V)$.

Recall that 
for a subset $A\subset X$, we call any function  $m:X\to\R_{\geq 0}$ with zero set $Z(m)=A$
a {defining function} for $A$.
We define a dg category $Open(X)$ as follows. The objects of $Open(X)$ are
pairs $\fU=(U,m)$ where $U\subset X$ is an open set, and $m:X\to \R_{\geq 0}$ is a defining function
for the complement $X\setminus U$.\footnote{In some contexts, 
it might be more natural to assume that
$m$ is a defining function for the boundary $\del U\subset X$. Our requirement that
$m$ vanishes on all of $X\setminus U$ plays no role outside of
reducing future notation.}
The complex of morphisms from an object $\fU_0=(U_0,m_0)$ to an object $\fU_1=(U_1,m_1)$ 
is the relative deRham complex
$$
\hom_{Open(X)}(\fU_0,\fU_1)=(\Omega(\ol U_0\cap U_1, \del U_0\cap U_1),d).
$$
Note the obvious fact that the morphisms are independent of the
defining functions. Given a third object $\fU_2=(U_2,m_2)$, the composition of
morphisms is the wedge product of forms
$$
\Omega(\ol U_0\cap U_1, \del U_0\cap U_1)
\otimes
\Omega(\ol U_1\cap U_2, \del U_1\cap U_2)
\to
\Omega(\ol U_0\cap U_2, \del U_0\cap U_2)
$$
To see this is well-defined, note that the support of any such wedge product
lies in $\ol U_1$, and thus since 
$$
(\del U_0\cap U_2) \cap \ol U_1
= (\del U_0\cap U_2 \cap U_1) \cup
(\del U_0\cap U_2 \cap \del U_1)
\subset (\del U_0\cap U_1) \cup
(U_2 \cap \del U_1),
$$
the support is disjoint from $\del U_0 \cap U_2$.

For an open subset $i:U\hookrightarrow X$, recall that 
$i_*\C_U$ denotes the standard extension of
the constant sheaf on $U$.

\begin{lem} \label{hom calc lemma}
For open subsets $i_0:U_0\hookrightarrow X,$ $i_1:U_1\hookrightarrow X$,
we have a canonical quasi-isomorphism
$$
\hom_{Sh(X)}(i_{0*}\C_{U_0},i_{1*}\C_{U_1})
\simeq
(\Omega(\ol U_0\cap U_1, \del U_0\cap U_1),d).
$$
The composition of morphisms coincides with the wedge product of 
differential forms.
\end{lem}

\begin{proof}
By standard identities, we have canonical quasi-isomorphisms
$$
\begin{array}{ccl}
\hom_{Sh(X)}(i_{0*}\C_{U_0},i_{1*}\C_{U_1}) 
& \simeq &
\Gamma(X, \intHom(i_{0*}\C_{U_0},i_{1*}\C_{U_1})) \\
& \simeq &
\Gamma(X, i_{1*}\intHom(i_1^*i_{0*}\C_{U_0},\C_{U_1}))) \\
& \simeq &
\Gamma(X, i_{1*}\intHom(\omega_{U_1},i_1^!i_{0!}\omega_{U_0})) \\
& \simeq &
\Gamma(X, i_{1*}i_1^*i_{0!}\C_{U_0}).
\end{array}
$$
Here we have written $\omega_{U_0},\omega_{U_1}$ for the dualizing complexes.
By de Rham's theorem, we also have a canonical quasi-isomorphism
$$\Gamma(X, i_{1*}i_1^*i_{0!}\C_{U_0})\simeq
(\Omega(\ol U_0\cap U_1, \del U_0\cap U_1),d).
$$

We leave it to the reader to check the last assertion.
\end{proof}

By the preceding lemma,
we may define a dg functor $P:Open(X)\to Sh(X)$ by sending an object $\fU=(U,m)$ to
the standard sheaf $i_*\C_U$
where $i:U\hookrightarrow X$ is the inclusion. 
By the preceding lemma and Proposition~\ref{generators}, 
the induced dg
functor on twisted complexes $\tw P:\tw Open(X)\to Sh(X)$ is a quasi-equivalence.


\subsection{Smooth boundaries}\label{smooth}

In what follows, we explain how to calculate morphisms in $Open(X)$
using open sets with smooth transverse boundaries.
To do this,
we will need to make choices of perturbation data.
It will be clear that the choices range over a contractible set,
and that they can be made compatibly
for any finite collection of objects.

Recall that the complex of morphisms from an object $\fU_0=(U_0,m_0)$ to an object $\fU_1=(U_1,m_1)$ 
is the relative deRham complex
$$
\hom_{Open(X)}(\fU_0,\fU_1)=(\Omega(\ol U_0\cap U_1, \del U_0\cap U_1),d).
$$
Our reinterpretation of this will be a complex not only quasi-isomorphic to it
but in fact isomorphic to it.

First, fix a Whitney stratification $\CS_0$ of $X$ compatible
with the boundary $\del U_0\subset X$,
and let $\Lambda_{\CS_0}\subset T^*X$ be the conical
conormal set
associated to $\CS_0$. By Lemma~\ref{sard},
there is $\ol \eta_1>0$ such that
there are no $\Lambda_0$-critical values of $m_1$ in the open interval $(0, \ol \eta_1)$.

\begin{lem}\label{first app of isotopy}
For $\eta_1\in(0,\ol\eta_1),$ there is a compatible collection of identifications
$$
(\ol U_0\cap X_{m_1>\eta_1}, \del U_0\cap X_{m_1 > \eta_1})
\simeq
(\ol U_0\cap U_1, \del U_0\cap U_1)
$$
which are the identity on $\ol U_0\cap X_{m_1\geq \ol\eta_1}$.
\end{lem}

\begin{proof}
By construction, there are no $\Lambda_{\CS_0}$-critical points of the map
$$
m_1: X_{0<m_1 <\ol\eta_1}\to (0,\ol\eta_1).
$$
For $\eta_1\in(0,\ol\eta_1),$ we may construct a compatible collection of diffeomorphisms
$$
(0,\ol\eta_1)\to (\eta_1,\ol\eta_1)
$$
by integrating an appropriate collection of vector fields.
Thus by the Thom isotopy lemma, we may lift these diffeomorphisms
to obtain identifications
$$
\ol U_0\cap X_{m_1>\eta_1}
\simeq
\ol U_0\cap U_1.
$$
Since $\CS_0$ is compatible with $\del U_0$,
the constructed identifications respect the pairs.
\end{proof}

Next, choose $\eta_1\in(0,\ol\eta_1)$,
fix the Whitney stratification $\CS_{\eta_1}$ of $X$ given by the hypersurface
$X_{m_1=\eta_1}$ and its complement,
and let $\Lambda_{\CS_{\eta_1}}\subset T^*X$ be the conical conormal set
associated to ${\CS_{\eta_1}}$. By Lemma~\ref{sard},
there is $\ol\eta_0>0$ such that
there are no $\Lambda_{\CS_{\eta_1}}$-critical values of $m_0$ in the open interval $(0,\ol\eta_0)$.

\begin{lem}\label{second app of isotopy}
For $\eta_0\in(0,\ol\eta_0),$ there is a compatible collection of identifications
$$
(X_{m_0>{\eta_0}}\cap X_{ m_1\geq \eta_1}, X_{m_0>{\eta_0}}\cap X_{m_1=\eta_1})
\simeq
(X_{m_0>{0}}\cap X_{ m_1\geq \eta_1}, X_{m_0>{0}}\cap X_{m_1=\eta_1}).
$$
which are the identity on $X_{m_0\geq \ol\eta_0}\cap X_{ m_1\geq \eta_1}$.
\end{lem}

\begin{proof} 
The argument is similar to that of the previous lemma.
By construction, there are no $\Lambda_{\CS_{\eta_1}}$-critical points of the map
$$
m_0: X_{0<m_0 <\ol\eta_0}\to (0,\ol\eta_0).
$$
For $\eta_0\in(0,\ol\eta_0),$ we may construct a compatible collection of diffeomorphisms
$$
(0,\ol\eta_0)\to (\eta_1,\ol\eta_0)
$$
by integrating an appropriate collection of vector fields.
Thus by the Thom isotopy lemma, we may lift these diffeomorphisms
to obtain identifications
$$
X_{m_0>{0}}\cap X_{ m_1\geq \eta_1}
\simeq
X_{m_0>{\eta_0}}\cap X_{ m_1\geq \eta_1}.
$$
Since $\CS_{\eta_1}$ is compatible with $X_{m_1=\eta_1}$,
the constructed identifications respect the pairs.
\end{proof}

Putting together the two previous lemmas, we obtain
a compatible collection of identifications
$$
U_0\cap U_1
\simeq
X_{m_0>{\eta_0}}\cap X_{ m_1> \eta_1}.
$$
Now we have the crucial observation:
for every open set $N_0\subset X$ containing $\del U_0\cap X_{m_1\geq \eta_1}$, 
there exists $\eta_0'>0$ such that $N_0$
contains the set
$X_{m_0<\eta_0'} \cap X_{ m_1\geq \eta_1}$.
Thus by construction, the above identifications induce a compatible collection of
isomorphisms of complexes
$$
(\Omega(\ol U_0\cap U_1, \del U_0\cap U_1),d)
\simeq
(\Omega(X_{m_0\geq \eta_0}\cap X_{m_1>\eta_1}, X_{m_0=\eta_0}\cap X_{m_1>\eta_1}),d).
$$
In what follows, we will use this reinterpretation of morphisms of $Open(X)$.


\subsection{Morse theory}\label{morse theory section}
In the previous discussion, we have constructed a dg functor
$$
P:Open(X) \stackrel{}{\to} Sh(X)
$$
such that $P$ identified $Sh(X)$ as a triangulated
envelope of $Open(X)$.
In this section, using Morse theory, we define an $A_\infty$-category $Mor(X)$ and an
$A_\infty$-functor
$$
M:Open(X)\stackrel{}{\to} Mor(X)
$$ 
which is a quasi-equivalence.

In this section, to simplify the exposition, we assume here that $X$ is {\em oriented},
and leave the general case to the reader.


\subsubsection{Manifolds with corners}
We begin by recalling some standard material from Morse theory.
We first discuss some general facts for an arbitrary open subset $W\subset X$, 
then specialize to the case where the closure $\ol W\subset X$ is a manifold
with corners. 

Let $W\subset X$ be an open subset.
Let $f:W\to \R$ be a function which extends to a small neighborhood of the closure
$\ol W\subset X$ such that
all critical points of the extension are nondegenerate and lie in $W$.
Let $Cr(f)\subset W$ denote the set of critical points, and let $i(x)$ denote the index of $x\in Cr(f)$.
Our convention is that a local minimum has index $0$, and a local maximum has
index $\dim X$.

Fix a Riemannian metric $g$ on $X$, and let $\nabla f$ denote the gradient vector field.
Let $\tilde W\subset W\times\R$ be a maximal domain
for the the gradient flow $\psi_t:\tilde W\to W$. For each $w\in W$, the fiber of $\tilde W$
above $w$
is an open (possibly unbounded) interval.
For each $x\in Cr(f)$, define the
stable and unstable manifolds
$$
W^-_x=\{w\in W| \lim_{t\to+\infty} \psi_t(w)=x\}
\qquad
W^+_x=\{w\in W| \lim_{t\to-\infty} \psi_t(w)=x\}.
$$
Implicit in the definition is that for $w\in W$ to lie in a stable or unstable manifold,
the fiber of $\tilde W$ above $w$ contains the appropriate half-line.
The stable manifold $W^-_x$ and unstable manifold $W^+_x$ are
diffeomorphic to balls of dimension $i(x)$ and $\dim X-i(x)$ respectively.
An orientation $\omega^-_x$ for the stable manifold $W^-_x$
 and an orientation $\omega^+_x$ for the unstable manifold $W^+_x$ are said
 to be compatible if at $x$
 the composite orientation
 $\omega_x^-\wedge \omega_x^+$ coincides with the ambient orientation of $X$.

Now we specialize to the case when the closure $\ol W\subset X$ is a manifold with corners.
To be precise, consider the quadrant
$$
\ol Q=\{(x_1,\ldots,x_n)\in \R^n| x_1\geq 0, x_2\geq 0\}.
$$
We assume that 
for each $w\in \ol W\subset X$, there is an open neighborhood
$N(w)\subset X$, an open set $U\subset \R^n$, and a $C^1$ diffeomorphism $\psi:N(w)\to U$
such that 
$$
\psi(\ol W\cap N(w))= \ol Q\cap U.
$$ 
Furthermore, we assume that
there are two smooth transverse
hypersurfaces $H_0,H_1\subset X$
such that $\del W\subset H_0\cup H_1$.

We will need the following notion of when a function $f$ on a manifold with
corners $\ol W\subset X$ and a Riemannian metric $g$ on $X$
are compatible.

\begin{dfn}
A pair $(f,g)$ where $f$ is a function on a neighborhood of $\ol W$, and $g$ is a Riemannian metric on $X$ is said to be {\em directed} 
if $(f,g)$ is Morse-Smale, and the gradient vector field $\nabla f$
is inward pointing along $H_0$ and outward pointing along $H_1$.
\end{dfn}

With the above set-up, if we have a directed pair $(f,g)$, 
then for each $x\in Cr(f)$,
the closures of the stable and
unstable manifolds satisfy
$$
\ol W^-_x\cap H_1=\emptyset
\qquad
\ol W^+_x\cap H_0=\emptyset.
$$


\subsubsection{Morse moduli spaces}
We next recall the moduli space of {\em gradient trees} from Fukaya-Oh~\cite{FO}. 

A {\em based metric ribbon tree} is a quadruple $(T,i,v_0,\lambda)$ of the following data. 
First, $T$ is a finite tree with $d+1$ end vertices and no vertex containing exactly
two edges. Second, $i:T\to D\subset \R^2$
is an embedding of $T$ in the closed unit disk 
such that the $d+1$ end vertices are precisely the intersection $i(T)\cap \del D$.
Third, $v_0$ is a distinguished end vertex of $T$. We refer to $v_0$ as the root vertex, and the other $d$ end vertices
as the leaf vertices. 
An edge $e\subset T$ is called an interior edge if $e$ does not contain an end vertex,
otherwise $e$ is called an exterior edge.
Finally, $\lambda$ is a function which assigns to every interior edge $e_{in}\subset T$ a positive
length $\lambda(e_{in})\in\R_+$. 
Two based metric ribbon trees are to be considered equivalent if there is an isotopy
of the closed disk which identifies all of the data.

Fix the orientation of the edges of $T$ such that all arrows point in
the direction of minimal paths from the leaf vertices to the root vertex.
The left and right sides of an edge refer to the components of the complement $D\setminus T$
which are respectively to the left and right of the edge with respect to the orientation of the edge.
Label the $d+1$ components of the complement $D\setminus T$ with elements of $\Z/(d+1)\Z$
in counterclockwise order
starting with $0$ for the component to the left of the edge terminating at the root vertex $v_0$.
For an edge $e\subset T$, let $\ell(e)$ and $r(e)$ denote the labelings of the
left and right sides of $e$.

For $i\in \Z/(d+1)\Z$, let
$U_i\subset X$ be an open subset with boundary 
$\del U_i\subset X$ a smooth hypersurface. 
Suppose that the boundaries $\del U_i$ form a transverse collection of hypersurfaces.
Let $f_i:U_i\to \R$ be a function which extends to a small neighborhood of the closure
$\ol U_i\subset X$. The difference $f_{i+1}-f_i$ is a function on the intersection
$U_i\cap U_{i+1}$ which extends to a small neighborhood of the closure
$\ol {U_i\cap U}_{i+1}\subset X$. 
Suppose the critical points of the extension of $f_{i+1}-f_i$
are nondegenerate and lie in $U_i\cap U_{i+1}$.
Suppose that we have chosen a Riemannian metric $g$ on $X$ such that
each pair $(f_{i+1}-f_i,g)$ is directed. In other words,
each gradient vector field $\nabla f_{i+1}-\nabla f_i$ points {inward}
along the hypersurface $\del U_{i+1}$ and {outward} along the
hypersurface $\del U_{i}$. 

With this set-up, a {\em gradient tree} is a pair $((T,i,v_0,\lambda),\tau)$ consisting
of a metric ribbon tree $(T,i,v_0,\lambda)$, and a continuous map $\tau:T\to X$ such that
the following holds.
\begin{enumerate}
\item For each end vertex $v\in T$, and the unique exterior edge $e_{ext}\in T$ containing it, we have
$$
\tau(v)\in Cr(U_{r(e_{ext})}\cap U_{\ell(e_{ext})},f_{r(e_{ext})}-f_{\ell(e_{ext})}).
$$
\item For each interior edge $e_{in}\subset T$, after making the identification 
$e_{in}\simeq[0,\lambda(e_{in})]$,
we have
$$
\tau'|_{e_{in}}=-\nabla (f_{\ell(e_{in})}- f_{r(e_{in})}).
$$
\item For each exterior edge $e_{ext}\subset T$, after making the identification 
$e_{ext}\simeq(-\infty,0]$,
we have
$$
\tau'|_{e_{ext}}=-\nabla (f_{\ell(e_{ext})}- f_{r(e_{ext})}).
$$

\end{enumerate}
Note that a single gradient tree alone contains the information of the images of the vertices,
and the oriented gradient lines which are the images of the edges. 

For each $i\in \Z/(d+1)\Z$, fix a critical point $x_{i}\in Cr(U_i\cap U_{i+1},f_{i+1}-f_i)$.
After a small perturbation of the data, the moduli space of gradient trees
$$
 \CM(T;f_{0},\ldots,f_{d};x_{0},\ldots,x_{d})
$$
with specified critical points is a manifold of dimension
$$
\sum_{i\in \Z/(d+1)\Z} i(x_i)- d \dim X + d-2.
$$
Orientations of the stable manifolds
of the differences $f_{i+1}-f_i$ induce a canonical orientation of the moduli space.
For example, for $d=1$, a single edge is the only based ribbon metric tree,
and the moduli space 
$
 \CM(T;f_{0},f_{1};x_{0},x_{1})
$
is the space of trajectories from $x_0$ to $x_1$
with orientations.


\subsubsection{Morse $A_\infty$-category}
Following Fukaya-Oh~\cite{FO}, we define an $A_\infty$-category $Mor(X)$ as follows.
As with $Open(X)$, 
the objects of $Mor(X)$ are
pairs $\fU=(U,m)$ where $U\subset X$ is an open set, and $m:X\to \R_{\geq 0}$ is a defining function
for the complement $X\setminus U$. To this data, we associate
the function $f:U\to \R$ defined by $f=\log m$.

To define the morphisms from an object $\fU_0=(U_0,m_0)$ 
to an object $\fU_1=(U_1,m_1)$, 
we will associate to them a directed pair and assign its Morse complex.
To ensure that we may find a directed pair,
we must refine the procedure summarized in Section~\ref{smooth}.

First, fix a Whitney stratification $\CS_0$ of $X$ compatible
with the boundary $\del U_0\subset X$,
and let $\Lambda_{\CS_0}\subset T^*X$ be the conical
conormal set
associated to $\CS_0$. By Lemma~\ref{sard},
there is $\ol \eta_1>0$ such that
there are no $\Lambda_0$-critical values of $m_1$ in the open interval $(0, \ol \eta_1)$.

Next, choose $\eta_1\in(0,\ol\eta_1)$,
fix the Whitney stratification $\CS_{\eta_1}$ of $X$ given by the hypersurface
$X_{m_1=\eta_1}$ and its complement,
and let $\Lambda_{\CS_{\eta_1}}\subset T^*X$ be the conical conormal set
associated to ${\CS_{\eta_1}}$. By Lemma~\ref{sard},
there is $\ol\eta_0>0$ such that
there are no $\Lambda_{\CS_{\eta_1}}$-critical values of $m_0$ in the open interval $(0,\ol\eta_0)$.

Now consider the critical points of the function $m_0\times m_1:X\to\R^2$. By definition,
critical points of $m_0\times m_1$ are points $x\in X$ such that $dm_0(x)$ and $dm_1(x)$
are colinear. Note that for points of $U_0\cap U_1$,
this is equivalent to $df_0$ and $df_1$ being colinear.
By construction, there are no critical values of $m_0\times m_1$
in the interval $(0,\ol\eta_0)\times \{\eta_1\}\subset \R^2$. In other words,
all critical points of $m_0\times m_1$ which lie on the hypersurface $X_{m_1=\eta_1}$
lie in the compact region $X_{m_0\geq \ol \eta_0,m_1=\eta_1}$.
Thus we may choose $\ol \epsilon_0>0$ such that for any $\epsilon_0\in (0,\ol \epsilon_0)$,
and any point $x\in X_{m_1=\eta_1}$ where $dm_0$ and $dm_1$ are colinear,
we have 
$$
\epsilon_0 \left |\frac{df_0(x)}{df_1(x)}\right |< 1.
$$
Here the fraction notation reflects the fact that the two covectors are colinear
so differ by a scalar.
Note in particular that $df_1(x)\neq 0$ on $X_{m_1=\eta_1}.$
This bound will guarantee that we may choose a Riemannian metric on $X$ such that
the gradient 
$\nabla f_1-\epsilon_0\nabla f_0$ is inward pointing along $X_{m_1=\eta_1}$.

Next, we need to choose $\eta_0>0$ small enough so that 
we may choose a Riemannian metric on $X$ such that
the gradient 
$\nabla f_1-\epsilon_0\nabla f_0$ is outward pointing along $X_{m_0=\eta_0}$.
Note that the values of $df_1$ along the compact hypersurfaces $X_{m_1=\eta_1}$
are bounded. Furthermore, the values of $df_0$ along the compact hypersurface
$X_{m_0=\eta_0}$ tend uniformly to infinity as $\eta_0$ tends to zero.
Thus for any $\eta_0>0$ small enough, 
and any point $x\in X_{m_0=\eta_0}$ where $dm_0$ and $dm_1$ are colinear,
we have 
$$
\left |\frac{df_1(x)}{\epsilon_0 df_0(x)}\right |< 1.
$$

In conclusion, for a sequence of sufficiently small choices $\eta_1>0$, then $\epsilon_0>0$,
then $\eta_0>0$, we have the following result.

\begin{lem}
There is a Riemannian metric $g$ on $X$
such that $(f_1-\epsilon_0 f_0,g)$ is a directed pair
on the manifold with corners $X_{m_0\geq \eta_0,m_1\geq \eta_1}$.
The set of such metrics is open and convex.
\end{lem}

\begin{proof}
The construction of the metric can be done locally using the bounds
of the above procedure.  The conditions on the metric are open and convex.
\end{proof}

Finally, we choose small perturbations of our functions and metric,
and define the space of morphisms to be the graded vector space generated
by critical points
$$
\hom_{Mor(X)}(\fU_0,\fU_1)=\C\{Cr(X_{m_0>\eta_0}\cap X_{m_1>\eta_1}, f_{1}-\epsilon_0 f_0)\}.
$$
The differential counts the oriented number of points 
in the moduli spaces of gradient trees
$$
m^1_{Mor(X)}(x_{0})=\sum_{T} \sum_{x_{1}} 
\deg \CM(T;\epsilon_0 f_{0},f_{1};x_{0},x_{1})\cdot x_{1}.
$$
with $T$ the interval $[-1,1]$.

\begin{lem}\label{morsecomplex}
$m^1_{Mor(X)}$ is well-defined and satisfies $(m^1_{Mor(X)})^2=0$.
\end{lem}

\begin{proof}
Fix critical points $x_0,x_1$.
Both assertions are implied by the claim that any sequence of gradient trajectories
$$\tau_n:[-1,1]\to X_{m_0> \eta_0,m_1> \eta_1}
$$ in the moduli space 
$\CM([-1,1];\epsilon_0 f_{0},f_{1};x_{0},x_{1})$ has a subsequence
which converges to a map 
$$\tau_\infty:[-1,1] \to X_{m_0> \eta_0,m_1> \eta_1}.
$$
The claim
with $\deg{x_1}=\deg{x_0}+1$ proves the first assertion. 
The claim
with $\deg{x_1}=\deg{x_0}+2$ implies that the boundary of the moduli space
in this case
is precisely the usual moduli space of broken trajectories which calculates $(m^1_{Mor(X)})^2$.
This  proves the second assertion.

To prove the claim, observe that the image set $\tau_n([-1,1])$ can never approach
the boundary 
$$
\del X_{m_0\geq  \eta_0,m_1\geq \eta_1} = 
(X_{m_0 = \eta_0} \cap X_{m_1 \geq \eta_1}) \cup
(X_{m_1 \geq  \eta_0} \cap X_{m_1 = \eta_1}).
$$
This follows form our assumptions on the behavior of the gradient vector field of $f_1-\epsilon f_0$
along the boundary (inward and outward pointing).
\end{proof}

By construction, the morphism complex $(\hom_{Mor(X)}(\fU_0,\fU_1),\mu^1_{Mor(X)})$
calculates the relative
cohomology 
$$
H^*(X_{m_0\geq \eta_0}\cap X_{m_1>\eta_1}, X_{m_0=\eta_0}\cap X_{m_1>\eta_1}).
$$

To define the higher compositions, for a finite collection of objects,
 we must follow the above procedure
sequentially. What we need is summarized in the following definition.

\begin{dfn}
\label{transversedef}
Consider a collection of pairs $(U_i,f_i)$ where $U_i\subset X$ is an open subset with boundary
$\del U_i$ a smooth hypersurface, and $f_i:U_i\to \R$ is a function indexed
by $i\in\Z/(d+1)\Z$. The collection said to be 
{\em transverse} if there is a Riemannian metric $g$ on $X$ such that the following holds.
For $i\in\Z/(d+1)\Z$, the hypersurfaces $\del U_i$ and $\del U_{i+1}$ are transverse,
and $(f_{i+1}-f_i,g)$ is a directed pair on $U_i\cap U_{i+1}$.
\end{dfn}

 In Section~\ref{a-infinity}, we will carefully explain in the context of the Fukaya category of $T^*X$
how to arrive at such a collection.
The procedure described there is modestly more complicated,
but strictly contains what is needed here.
Therefore we will not pursue further details here, but only mention
the following salient points.

Given a collection of objects indexed by  $i\in \Z/(d+1)\Z$, for any 
sequence of sufficiently small choices $\epsilon_i>0$ and $\eta_i>0$,
one may arrange for the perturbed boundaries $X_{m_i=\eta_i}$ to form a transverse
collection.
Furthermore, one may sequentially obtain bounds on the differentials $df_i$ along
the perturbed boundarues.
Together this allows one
to find dilations and a Riemannian metric on $X$ such that all dilated pairs are directed.
After performing small perturbations,
the higher composition maps count the oriented number of points 
in the moduli spaces of gradient trees
$$
m^d_{Mor(X)}(x_{0},\ldots,x_{d-1})=\sum_{T} \sum_{x_{d}} 
\deg \CM(T;\epsilon_0 f_{0},\ldots,\epsilon_{d-1} f_{d-1},f_{d};x_{0},\ldots,x_{d})\cdot x_{d}.
$$
In the following section, we will apply homological perturbation theory to verify
the following.

\begin{prop}\label{a-infinity morse relations}
The maps $\mu^d_{Mor(X)}$ are well-defined and satisfy the $A_\infty$-quadratic composition rule.
\end{prop}

In conclusion, it is worth commenting about the choices involved in the construction
of $Mor(X)$.
For a collection of objects indexed by $i\in \Z/(d+1)\Z$,
there are the small choices of constants $\eta_i>0$ to obtain smooth boundaries
and $\epsilon_i>0$ to dilate functions. These may be organized into a contractible ``fringed set"
as discussed in  Section~\ref{lagsec}. 
In addition, there is the choice of a Riemannian metric 
to obtain directed pairs. While not a perturbation in any sense, such metrics form a convex set.
Finally, there are the small perturbations of the functions and metric. This is no
different from the standard context.


\subsubsection{From differential forms to Morse theory}\label{homol pert theory applied}
Following Kontsevich-Soibelman~\cite{kontsoib},
we apply here the formalism of homological perturbation theory to prove
Proposition~\ref{a-infinity morse relations} and obtain
an $A_\infty$-functor
$$
M:Open(X)\stackrel{}{\to} Mor(X)
$$
which is a quasi-equivalence.
We will apply the formalism in the special case
of an idempotent.
The construction of the idempotent and the homotopy 
is completely analogous to that of Kontsevich-Soibelman:
one composes the limit operators of Harvey-Lawson~\cite{HL} with a smoothing
operator. 
To explain this, we return to the general context of a submanifold with corners $\ol W\subset X$
and boundary hypersurfaces $H_0,H_1\subset \del W$
with which we began 
this section.

Let $D'(W,H_0)$ denote the space of currents dual to the space of differential forms $\Omega(W, H_1)$.
There are two simple ways to obtain elements of $D'(W,H_0)$. First, we have the inclusion
$$
i:\Omega(W,H_0)\to D'(W,H_0)
$$
defined by taking the wedge product of forms and integrating.
Second, any oriented closed submanifold $V\subset W$ satisfying $\ol V \cap H_0=\emptyset$ defines 
an element 
$$[V]\in D'(W,H_0)$$ 
by integration. In particular, for each $x\in Cr(f)$,
the unstable manifold $W^+_x$ with a given orientation defines a current 
$[W^+_x]\in D'(W,H_0).
$
Similarly, the stable manifold $W^-_x$ with a given orientation defines a current 
$[W^-_x]\in D'(W,H_1)$.

In this context, the main result of Harvey-Lawson takes the following form.
Define the linear operator 
$
p:\Omega(W, H_0)\to D'(W, H_0)
$
by the kernel
$$
\mathcal P=\sum_{x\in Cr(f)} [W^-_x]\times [W^+_x]
$$
where the stable and unstable manifolds are given compatible orientations.
Define the homotopy operator
$
h:\Omega(W,H_0)\to D'(W, H_0)
$
by the kernel
$$
\mathcal H=\bigcup_{0\leq t< +\infty} [\Gamma_{\psi_t}]
$$
where $\Gamma_{\psi_t}\subset X\times X$ denotes the graph of 
the gradient flow $\psi_t$.
Then we have the equation of operators
$$
p-i=
d h+hd. 
$$

Now to obtain an honest idempotent $\pi$ and homotopy operator $T$
on $\Omega(W,H_0)$, we need only
compose with a smoothing operator $D'(W,H_0)\to \Omega(W,H_0)$. The details
of this are no different from the case considered by Kontsevich-Soibelman.

Applying the formalism of homological perturbation theory --
and recognizing that it coincides with counting gradient trees -- 
we 
both see that Proposition~\ref{a-infinity morse relations} must hold and obtain
an $A_\infty$-functor
$$
M:Open(X)\stackrel{}{\to} Mor(X)
$$
which is a quasi-equivalence.


\section{The Fukaya category}

The Fukaya $A_\infty$-category $Fuk(M)$ of a symplectic manifold $M$ is
a quantization of the Lagrangian intersection theory of $M$.
Roughly speaking,
its objects are Lagrangian submanifolds and its morphisms are generated by intersection points
of the Lagrangians. Its composition maps are defined by choosing
a compatible almost complex structure,
and counting 
holomorphic
polygons 
with boundary lying on the Lagrangians.\footnote{Properly, we should say ``pseudoholomorphic,'' but
we omit the prefix throughout.}
For example, for intersection points $p_{0} \in
L_0\cap L_1$ and $p_{1} \in L_1\cap L_2$, the coefficient of $p_{2}\in L_2\cap L_0$
in the product $p_{0}\cdot p_{1}$ is the number of holomorphic maps to $M$
from a disk with three marked boundary points mapping to the intersection points and
with the arcs between them mapping to the Lagrangians.

To this coarse description there are many details, refinements, and specializations
for various settings.  In this paper, we will use a composite picture
of the treatments from Eliashberg-Givental-Hofer \cite{SFT},
Fukaya-Oh-Ohta-Ono \cite{FOOO}, and Seidel \cite{Seidel}.
Our symplectic manifold is the cotangent bundle $M=T^*X$ of a compact real analytic manifold $X$.
We equip $T^*X$ with the exact symplectic form $\omega = d\theta,$
where $\theta$ is the canonical one-form 
$$\theta(v)\vert_{(x,\xi)} = \xi(\pi_*v),$$
with $\pi:T^*X\to X$ the standard projection.
For any choice of Riemannian metric on $X$, 
the associated Levi-Civita connection provides a compatible almost complex structure on $T^*X$,
along with a canonical Riemannian (Sasaki) metric on $T^*X$.  
We will also
consider a mild variation of these structures as described in Section \ref{warp}
below.

In what follows, we focus on the aspects of our situation which deviate from what is by now standard
in the subject. All of these differences stem from the fact that we will allow  closed but
noncompact Lagrangians.

We often use the following notation:
given a space $Y$, a function $g:Y\to \R$, and $r\in \R$, we write $Y_{g=r}$ for the subset
$\{y\in Y| g(y)=r\}$, and similarly for inequalities.


\subsection{Basics of $T^*X$}
\label{basics}

\subsubsection{Compactification}

Consider the bundle 
$J^1_{\geq 0}(X)=T^* X \times \R_{\geq 0}
$
of $1$-jets of non-negative
functions on $X$, and
let 
$J^1_{\geq 0}(X)'=J^1_{\geq 0}(X)\setminus (X \times \{0\})$
be the complement of the zero section. 
The multiplicative group $\R_{+}$ acts freely on $J^1_{\geq 0}(X)'$
by dilations.
The quotient
$$
\ol T^* X = J^1_{\geq 0}(X)'/\R_+
$$
equipped with the obvious projection $\ol\pi:\ol T^*X\to X$
provides a relative compactification of $\pi:T^*X\to X$.
We have the canonical inclusion $T^*X\hookrightarrow \ol T^* X$ which sends
a covector $\xi$ to the class of $1$-jets $[\xi,1]$, and we refer
to this inclusion implicitly whenever we consider $T^*X$ as a subset of $\ol T^*X$.
The divisor at infinity 
$$T^\infty X=\ol T^*X\setminus T^*X
$$ 
consists of the class of
$1$-jets of the form $[\xi, 0]$ with $\xi$ a non-zero covector.

Let $\CO(-1)$ denote the tautological $\R_+$-principal bundle $J^1_{\geq 0}(X)'\to \ol T^* X$ with fiber $\R_+\cdot(\xi,r)$
at the point $[\xi,r]$. 
The restriction of $\CO(-1)$ to the open subset $T^*X$ is canonically trivialized
by the section $[\xi,r]\mapsto (\xi/r,1)$.
Let $\CO(-1)_\infty$ denote the restriction of $\CO(-1)$ to the divisor at infinity $T^\infty X$.
A choice of trivialization of $\CO(-1)_\infty$ 
is equivalent to
a choice of (co-)sphere subbundle $S^*X\subset T^*X$, and provides a
canonical identification 
$$
T^\infty X\simeq S^*X.
$$
In fact, it is always possible to trivialize $\CO(-1)$ itself over all of $\ol T^*X$.
For example, if we choose a Riemannian metric on $X$, then 
we have the section 
$$[\xi,r]\mapsto (\hat \xi, \hat r),
\mbox{ where }|\hat\xi|^2 + \hat r^2=1.
$$
This identifies $\ol T^*X$ with the closed unit disk bundle,
and $T^\infty X$ with its unit sphere bundle. 
Note that such a trivialization can not be made equal to
the canonical trivialization of $\CO(-1)$
over the open subset $T^*X$.

By working with a spherical compactification
rather than a projective compactification, we pay the price of dealing with a manifold
with boundary. But we choose this approach because we will encounter objects on $T^\infty X$
which are not invariant under the antipodal involution.

The pull-back of $\theta$ to $J^1_{\geq 0}(X)$ descends to  
a one-form $\ol \theta$ on $\ol T^* X$
with values in the $\R$-line bundle $L(1)$ associated to the $\R_+$-principal bundle dual to $\CO(-1)$.
The restriction $\theta^\infty=\ol\theta|_{T^\infty X}$ is an $L(1)$-valued contact form on $T^\infty X$.
The canonical trivialization of $\CO(-1)$ over the open subset $T^* X$ identifies
the restriction $\ol\theta|_{T^*X}$ with the original one-form $\theta$.
By choosing a trivialization of $\CO(-1)_\infty$, we may consider $\theta^\infty$ as an honest
contact form. Equivalently, by choosing a sphere bundle $S^*X\subset T^*X$, we may identify
$\theta^\infty$ with the restriction of $\theta$ to $S^*X$. 
If we do not fix such identifications, we still have a well-defined contact 
structure $\ker(\theta^\infty)\subset T T^\infty X$, and a well-defined notion of positive normal direction.
This positive direction is an example of a structure on $T^\infty X$ which is not invariant
under the antipodal map.


\subsubsection{Geodesic flow}

Given a function $H:T^*X\to \R$, we have the Hamiltonian vector field $v_H$ defined
by 
$$dH(v)=\omega(v,v_H).
$$
When possible, integrating $v_H$ provides a Hamiltonian isotopy $\varphi_{H,t}:T^*X\to T^*X$.

A Riemannian metric on $X$ defines a Riemannian
(Sasaki) metric on $T^*X$, and thus an identification
$T^*(T^*X) \simeq T(T^*X)$. 
The canonical one-form $\theta$ on $T^*X$ corresponds to 
 the geodesic vector field $v_\theta$. 
  On the complement of the zero section $T^*X \setminus X$,
we have the normalized geodesic vector field
$\hat v_\theta= v_\theta/|v_\theta|$. It is the Hamiltonian vector field for the 
length function $H:T^*X \setminus X\to \R$ given by $H(x,\xi)=|\xi|$.
We write $\gamma_t:T^*X\setminus X\to T^*X\setminus X$ for the normalized geodesic flow for time $t$ associated to $\hat v_\theta$.
By definition, if we identify a covector $(x,\xi)\in T^*X$ with a vector $(x,v)\in TX$,
then we have the identity
$$
\gamma_t(x,v)=\exp_{x,t}(\hat v)_* (v)
$$
where $\hat v = v/|v|,$
the map $\exp_{x,t}:T_{x} X\to X$ denotes the exponential flow from the
point $x$ for time $t$, and the asterisk subscript indicates the
derivative (push-forward).
Since $v_\theta$ grows at infinity, its flow does not have a well-defined limit.
But $\hat v_\theta$ extends to give a Reeb
flow on the contact manifold at infinity $T^\infty X.$

A function $H:T^*X\to \R$ is said to be {\em controlled} if there is a compact set $K\subset T^*X$ such that outside of $K$ we have $H(x,\xi)=|\xi|$. 
The corresponding Hamiltonian isotopy 
$\varphi_{H,t}:T^*X\to T^*X$ equals the normalized geodesic flow $\gamma_t$
outside of $K$. 
Note that 
for any controlled function $H$,
the vector field $v_H$ may be integrated to a Hamiltonian isotopy $\varphi_{H,t}$, for all times $t.$  
Note as well that Hamiltonian flow by $|\xi|$ depends on the metric on $X$, but is 
independent of any choice of metric on $T^*X$.


\subsubsection{Almost complex structures}
\label{warp}
To better control holomorphic disks in $T^*X$, it is useful to introduce 
an almost complex structure $J_{con}$ which  near infinity is invariant under dilations.

Recall that a Riemannian metric on $X$ provides a canonical splitting
$$
T(T^*X) \simeq T_{b} \oplus T_{f},
$$
where $T_b$ denotes the horizontal base directions and $T_f$ the 
vertical fiber directions,
along with a canonical isomorphism 
$
j_0:T_b \stackrel{\sim}{\to} T_f
$
of vector bundles.
Thus we can define a compatible almost complex structure $J_{Sas}$ by the matrix
$$
J_{Sas}=
\left(
\begin{matrix}
0 & j_0^{-1} \\
-j_0 & 0
\end{matrix}
\right).
$$
We refer to $J_{Sas}$ as the Sasaki almost complex structure,
since by construction, the Sasaki metric 
is given by $g_{Sas}(U,V) = \omega(U, J_{Sas} V)$.

Fix a positive function $w:T^*X \rightarrow \R_{>0}$, and
define a new compatible almost complex structure $J_w$ by the matrix
$$
J_w=
\left(
\begin{matrix}
0 & w^{-1} j_0^{-1} \\
-w j_0 & 0
\end{matrix}
\right).
$$
Fix a local orthonormal frame $\{b_i, f_i\}_{i=1}^{\dim X}$ 
for $T_b\oplus T_f$ with respect to the Sasaki metric.
Then with respect to the new metric $g_w(U,V) = \omega(U, J_w V)$,
the lengths of the Sasaki frame take the form
$$
|b_i|_{g_w} =
w^{1/2} 
\qquad
|f_i |_{g_w} =w^{-1/2}. 
$$

For concreteness, we specialize the construction
by making a specific choice of the function $w$.
Namely, fix positive constants $r_0, r_1>0$,
and a bump function $b: \R\to \R$ such that $b(r) = 0$ for $r <r_0$, and $b(r) = 1$, for $r>r_1$.
Fix a constant $\beta\in \R$,
and set 
$$w(x,\xi) = |\xi|^{\beta b(|\xi|)},
$$ where as usual
$|\xi|$ denotes the length of a covector with respect to the original
metric on $X$. 
In particular, when $\beta = 0$, we recover the original Sasaki almost complex structure
$J_{Sas}$ and Sasaki metric $g_{Sas}$.

In what follows, we will restrict our attention to the choice $\beta = 1$
and denote the almost complex structure by $J_{con}$ and corresponding metric by $g_{con}$.
We will refer to $J_{con}$ as the conical almost complex structure since
near infinity, $J_{con}$ is invariant under dilations, and
the lengths of our Sasaki frame take the form
$$
|b_i|_{g_{con}} =
|\xi|^{1/2} 
\qquad
|f_i |_{g_{con}} = |\xi|^{-1/2}.
$$
Thus near infinity, we have replaced the Sasaki geometry
with a cone over the unit {(co-)sphere} bundle $S^*X$. To be precise,
if $ds^2$ is the restriction of the Sasaki metric to $S^*X$,
and for clarity we write $r$ for the length $|\xi|$,
then near infinity we have
$$
g_{con} = r^{-1/2} dr^2 + r^{1/2} ds^2.
$$
(Substituting $\tilde r=r^{1/2}$, one sees the familiar presentation of the metric of a cone.)
It is straightforward to check that $g_{con}$ is complete, and the normalized geodesic flow $\gamma_t$
is an isometry.

We will confirm 
in Section~\ref{hol disks} that 
the almost complex structure $J_{con}$ provides compact moduli spaces of holomorphic disks in the circumstances under consideration.
One can view the metric $g_{con}$ as being compatible with the compactification $\ol T^*X$
in the sense that near infinity, 
it treats base and angular fiber directions on equal footing: near infinity, the metrics
on the level sets of $|\xi|$ are simply scaled by the factor $|\xi|^{1/2}$.


\subsection{Lagrangians}
\label{lagsec}
Fix an analytic-geometric category $\CC$. 

\begin{lem}\label{disccritvals}
For any $\CC$-subset $\ol V\subset \ol T^*X$, there exists $r>0$ such that $|\xi|$ has no
critical points on $\ol V\cap T^* X$ for $|\xi|\geq r.$
\end{lem}
\begin{proof} The critical values of $1/|\xi|$ are a discrete $\CC$-subset of $\R$.
\end{proof}

\begin{lem}\label{flowdown}
Let $W$ be a compact space, and let $V\subset T^*X$ be a subset such that $|\xi|$
has no critical points on $V$ for $|\xi|\geq r$. Then any map $W\to V$ is homotopic in $V$
to a map $W\to V_{|\xi|<r}$
\end{lem}

\begin{proof} By the Thom
isotopy lemma,
we may use the gradient of $|\xi|$ to flow the image of the map $W\to V$.
\end{proof}

A subset $V\subset T^*X$ is said to be {\em conical}
 if it is invariant under the action of $\R_+$ by fiberwise dilations.

\begin{lem} 
If $V$ is a conical $\omega$-isotropic subset of $T^*X$, then
$\ol V\cap T^\infty X$
is a $\theta^\infty$-isotropic subset of $T^\infty X$.
\end{lem}

\begin{proof}
The one-form $\theta$ may be obtained from the symplectic
form $\omega$ by contracting with the Liouville vector field $v_\theta$.
The action of $\R_+$ by dilations is generated by $v_\theta$.
\end{proof}

As long as we assume that $\ol V\subset \ol T^*X$ is a $\CC$-subset,
we have the following very general result.

\begin{lem} 
If $\ol V\subset \ol T^*X$ is a $\CC$-subset such that $\ol V\cap T^*X$ is an $\omega$-isotropic
subset of $T^*X$,
then $\ol V\cap T^\infty X$
is a $\theta^\infty$-isotropic subset of $T^\infty X$.
\end{lem}

\begin{proof}
Let $\CN^\infty X$ be the family with general fiber $\ol T^*X$ and special fiber
the normal cone $N^\infty X$ of $\ol T^*X$ along the divisor $T^\infty X$. 
Let $\ol C\subset N^\infty X$ be the limit of $\ol V$ under specialization in the family $\CN^\infty X$.
(See \cite{KS}, pp. 185--187, for an exposition of the normal cone
and the specialization of subsets.)
By construction, $\ol C$
is a conical subset satisfying
$\ol C \cap T^\infty X = \ol V \cap T^\infty X.$

We claim that $\ol C$ is $\omega$-isotropic.
(The normal cone $N^\infty X$ inherits a well-defined $\omega$-isotropic
distribution.)
If we can show this, then we are done by the previous lemma.
To see this, choose a
Whitney stratification of $\CN^\infty X$ compatible with $\ol C$;
this is possible since $\ol V$ is a $\CC$-subset.
Then the Whitney condition and the fact that being $\omega$-isotropic is a closed
condition together imply the assertion:
the tangent spaces of the limit
$\ol C$ are contained
in the limits of the $\omega$-isotropic tangent spaces of $\ol V$.
\end{proof}

We will need to separate $\theta^\infty$-isotropic subsets of $T^\infty X$ using
the normalized geodesic flow (Reeb flow). 
To organize this, we use a variant of the notion of a {\em fringed set} from~\cite{GM}.
To define what a fringed set $R_{d+1}\subset \R^{d+1}_+$ is, we proceed inductively.
A fringed set $R_1\subset\R_+$ is any interval of the form $(0,r)$ for some $r>0$. 
A fringed set $R_{d+1}\subset\R_+^{d+1}$ is a subset satisfying the following:
\begin{enumerate}
\item $R_{d+1}$ is open in $\R^{d+1}_+$.
\item Under the projection $\pi:\R^{d+1}\to \R^d$, the image $\pi(R_{d+1})$ is a fringed set.
\item If $(r_1,\ldots, r_d, r_{d+1})\in R_{d+1}$, then $(r_1,\ldots, r_d, r'_{d+1})\in R_{d+1}$
for all $0<r'_{d+1}< r_{d+1}$.
\end{enumerate}
It is easy to check that fringed sets as defined here are contractible.

\begin{lem}
\label{deform}
For $i=0,\ldots, d$, let $V^\infty_i\subset T^\infty X$ be $\theta^\infty$-isotropic
compact subsets.
Then there is a fringed set $R_{d+1}\subset \R^{d+1}$ such that for $(\delta_0,\ldots,\delta_d)\in R_{d+1}$,
the normalized geodesic flow (Reeb flow) separates the subsets:
$$
\gamma_{\delta_i}(V^\infty_i)\cap 
\gamma_{\delta_j}(V^\infty_j) =\emptyset,
\mbox{ for $i\not =j$}.
$$
\end{lem}

\begin{proof}
We prove the assertion by induction. For $d=0$, there is nothing to prove.
Suppose we know the assertion for $d-1$ and seek to establish it for $d$.
For $(\delta_0,\ldots,\delta_{d-1})\in R_{d}$, consider
the $\theta^\infty$-isotropic subset
$$
V^\infty_{<d}=\bigcup_{i<d} \gamma_{\delta_i}(V_i^\infty).
$$
It suffices to show that there is $\delta_d=\delta_d(\delta_0,\ldots,\delta_{d-1})>0$
such that for all $0<\delta'_d<\delta_d$, we have
$$
\gamma_{\delta'_d}(V^\infty_d)\cap V^\infty_{<d}
=\emptyset.
$$
Suppose this were not true. Then by the curve selection lemma (Proposition \ref{csl}),
there is a $C^1$ map $\alpha:[0,1)\to V^\infty_d$ such that for all $t\in (0,1),$
$$
\gamma_t(\alpha(t))\in \gamma_t(V^\infty_d)\cap V^\infty_{<d}.
$$
In particular, $\gamma_t(\alpha(t))$ lies in the $\theta^\infty$-isotropic subset $V^\infty_{<d}$.
But we calculate
$$
\left.\frac{d}{dt}\gamma_t(\alpha(t))\right|_{t=0}=\gamma'_0(\alpha(0)) + (\gamma_0)_*\alpha'(0).
$$
Since $\alpha(t)$ lies in the $\theta^\infty$-isotropic subset $V^\infty_d$,
$\alpha'(t)$ is in the kernel of $\theta^\infty$ and
we arrive at the conclusion 
$$
\theta^\infty(\gamma'(\alpha(0)))=0.
$$
But this quantity is nonzero since $\gamma'$ is the Reeb vector field on $T^\infty X$.
\end{proof}


\subsubsection{Exact Lagrangians}

A Lagrangian $i : L\hookrightarrow T^*X$
is said to be {\em exact} if the restriction $i^*\theta$ is an
exact differential form.

\begin{lem}
\label{zeroarea}
Let $L\subset T^*X$ be a Lagrangian, $\Sigma$ a compact 
Riemann surface with boundary $\del\Sigma$,
and $u:(\Sigma,\del \Sigma) \to (T^*X,L)$ a differentiable map. Let
$A(u) = \int_\Sigma u^*\omega$ denote its symplectic area.
Then we have
\begin{enumerate}
\item $A(u)$
depends only on the homotopy class in $L$ of $u|_{\del\Sigma}$.
\item $A(u_\epsilon) = \epsilon A(u),$ where
$u_\epsilon: (\Sigma,\partial \Sigma)\to (T^*X,\epsilon L)$
is the composition of $u$ with the dilation $(x,\xi)\mapsto (x,\epsilon \xi).$
\item
If $L$ is exact, $u$ is constant.
\item If $u$ is holomorphic, then $Area(u) = A(u).$
\end{enumerate}
\end{lem}

\begin{proof}
For the first assertion, note that if $\gamma_1$ and $\gamma_2$
are homotopic loops in $L$ and $S\subset L$ satisfies $\partial S = [\gamma_1] - [\gamma_2],$
then $\oint_{\gamma_1}\theta-\oint_{\gamma_2}\theta
= \int_S d\theta = \int_S \omega = 0,$ since $S\subset L$ is $\omega$-isotropic.
To prove the second,
for $p = (x,\xi)\in T^*X,$
let $\epsilon p$ denote the point $(x,\epsilon \xi)$,
and note that $\theta\vert_{\epsilon p}(\epsilon v) = \epsilon \theta\vert_p v.$
The third claim follows from exactness:
$\int_\Sigma u^*\omega = \oint_{\partial \Sigma} u^*\theta = \oint_{\partial \Sigma}
u^* (d\psi) = 0.$  The fourth statement expresses the
fact that when $J$ is a compatible almost complex structure,
$J$-holomorphic maps are calibrations for $\omega.$
\end{proof}


\subsubsection{Standard Lagrangians}

Let $Y\subset X$ be a submanifold. The conormal bundle $T^*_Y X\subset T^*X$
is homotopic to its zero section $Y$, and thus
is an exact Lagrangian, since $\theta$ is identically zero on the zero section $X$.

Given a defining function $m:X\to\R_{\geq 0}$
for the boundary $\del Y\subset X$, we define
$f:X\setminus \partial Y \to \R$ by $f=\log m$ and define
the {\em standard Lagrangian} $L_{Y,f}\subset T^*X\vert_Y \subset T^*X$ to be
the fiberwise sum
$$
L_{Y,f}= T^*_Y X + \Gamma_{df}\vert_Y,
$$
where $\Gamma_{df}\subset T^*X\vert_{X\setminus \partial Y}$
denotes the graph of $df,$
and the sum is taken fiberwise in $T^*X\vert_Y.$
By construction,
$L_{Y,f}$ depends only on the restriction
of $m$ to $Y$: if two functions agree on $Y$, then over $Y$ their differentials differ by
a section of the conormal $T^*_Y X.$  For this reason, in the sequel we will
often refer to
$m$ and $f$ as functions on $Y.$
Note that if $Y$ is an open submanifold,
we could equivalently take $m$ to be a defining function
for the complement $X\setminus Y$.

\begin{lem}
$L_{Y,f}$ is canonically Hamiltonian isotopic to $T^*_Y X$. In particular,
$L_{Y,f}$ is exact.
\end{lem}

\begin{proof}
To avert potential confusion, it is worth pointing out that $T^*_Y X$ is not
necessarily closed in $T^*X$, and we will move $L_{Y,f}$ through a family
of Lagrangian submanifolds which are not necessarily closed.
In the subset $T^*X|_Y\subset T^*X$, consider the function $H=f\circ \pi$ and the associated Hamiltonian flow $\varphi_{H,t}$. 
One checks that $\varphi_{H,t}$ applied to the Lagrangian $L_{Y,f}$ takes it to its dilation $(1-t)\cdot L_{Y,f}$. In particular, when $t=1$,
one arrives at the conormal Lagrangian $T^*_Y X$.
\end{proof}


\subsection{Brane structures}
\label{branes}

In order to define a Fukaya category of a symplectic
manifold $M$, one needs 
a grading on the
Lagrangian intersections and 
orientations of the relevant
moduli spaces of holomorphic disks. 
(Alternatively, one could be satisfied with
an ungraded version of the Fukaya category with
characteristic $2$ coefficients.)
Topological obstructions to gradings come from
the bicanonical bundle of $M$
and the Maslov class of Lagrangians.
Orientation of the moduli spaces requires a relative pin structure
on the Lagrangians, so that their second Stiefel-Whitney
classes must be restrictions of a (common) class on $M.$\footnote{Some
authors have proven orientability under more restrictive conditions.
In \cite{Seidel}, Lagrangians are assumed to be pin,
while in \cite{FOOO} they are taken to be oriented and relatively spin.  Our
more general condition follows Wehrheim and Woodward's work in progress.
We note that the Lagrangians of interest for us will be canonically pin
when $X$ is pin, and canonically oriented and relatively spin when $X$
is oriented.}
In this section, we show all
obstructions to these structions vanish for $M=T^*X$ and the Lagrangians of
interest.
In what follows, we always work with the
canonical exact symplectic structure on $T^*X$,
and the compatible almost complex structure induced by a
Riemannian metric on $X$.


\subsubsection{Bicanonical line}
The almost complex structure on $T^*X$ allows us to define the holomorphic canonical bundle
$$
\kappa= (\wedge^{\dim X}T^{hol}T^*X)^{-1}.
$$
In order to compare the squared phase of Lagrangian subspaces at different points
of $T^*X$, we need a homotopy class of trivializations of the bicanonical bundle $\kappa^{\otimes 2}$.

\begin{prop} The bicanonical bundle $\kappa^{\otimes 2}$ of $T^*X$ is canonically trivial.
\end{prop}

\begin{proof}
Since the zero section $X$ is a deformation retract of $T^*X,$
it suffices to see $\kappa^{\otimes 2}|_X$ is canonically trivial.
At the zero section, $TT^*X$ has a canonical splitting into vertical and
horizontal spaces,
$$
TT^*X\vert_X = T^*X\oplus TX\vert_X \cong TX\otimes \C\vert_X,
$$
where we have identified the cotangent bundle with the tangent
bundle using the metric, and identified the normal directions with
the imaginary directions using the compatible almost complex structure.
As a result, 
$$T^{hol}(T^*X)\vert_X \cong TX\otimes \C,
$$ 
and we
see 
$$
\kappa\vert_X \simeq \pi^*(or_X)\otimes \C,$$
where $or_X$ is the orientation line bundle.  Thus
$\kappa$ is trivializable if and only if $X$ is orientable,
and $\kappa^{\otimes 2}$ is canonically trivial for any $X$.
\end{proof}

\begin{rmk}
In general, trivializations of a complex line bundle over a space $X$
form a torsor over the group of maps $X\to \C^*$. Homotopy classes of trivializations
form a torsor over the group $H^1(X,\Z)$. 
\end{rmk}


\subsubsection{Grading}

Let $\eta^2$ be the canonical trivialization of $\kappa^{\otimes 2}$,
and let ${\mathcal Lag}_{T^*X}\to T^*X$ be the bundle of Lagrangian
planes.
We have the squared phase map
$$
\alpha : {\mathcal Lag}_{T^*X} \rightarrow U(1)
$$
$$
\alpha(\mathcal L)= \eta(\wedge^{\dim X} \mathcal L)^2/ |\eta(\wedge^{\dim X} \mathcal L)|^2.
$$

For a Lagrangian $L\subset T^*X$ and a point $x\in L$, we obtain a map $\alpha:L\to U(1)$
by setting $\alpha(x)=\alpha(T_xL).$
The Maslov class $\mu(L)\in H^1(L)$
is the obstruction class
$$
\mu = \alpha^*(dt),
$$ 
where $dt$ is the standard one-form on $U(1)$.
Thus $\alpha$ has a lift to a map $\widetilde\alpha:L\to \R$
if and only if $\mu = 0.$
Such a lift is called a {\em grading} of the Lagrangian.

\begin{rmk}
Choices of gradings of a Lagrangian $L$  form a torsor over the group
$H^0(L,\Z)$. 
\end{rmk}

Next we check that our standard Lagrangians have canonical gradings.  
Recall that to a submanifold $Y\subset X$, and a defining function $m:X\to\R_{\geq 0}$
for the boundary $ \del Y\subset X$, we have the standard Lagrangian 
$$
L_{Y,f}= T^*_Y X + \Gamma_{df} \subset T^*X
$$ 
where $f:Y\to \R$ is given by $f=\log m$.

\begin{prop}
The Maslov class $\mu(L_{Y,f})\in H^1(L_{Y,f})$ vanishes. In fact, there is 
a canonical grading of $L_{Y,f}$.
\end{prop}

\begin{proof}
Since $L_{Y,f}$ is canonically Hamiltonian isotopic to $T^*_Y X$, it suffices to check the
assertions for $T^*_Y X$. Furthermore, since $T^*_Y X$ is a vector bundle
over $Y$, it suffices to check the assertions along $Y$.
Let $\{e_j\}_{j=1}^{\dim X}$ be an orthonormal frame field for $X$ along $Y$
extending
an orthonormal frame field $\{e_j\}_{j=1}^{\dim Y}$ for $Y$. 
Note that the zero section, and in particular the frame field $\{e_j\}_{j=1}^{\dim X}$,
has constant squared phase equal to the identity of $U(1)$.
Thus we can equip it with the canonical constant grading given by $0$ in $\R$.
Then the frame field along $Y$ given by 
$
\{e_j\}_{j=1}^{\dim Y} \cup \{Je_j\}_{j=\dim Y+1}^n
$
has constant squared phase $(-1)^{\codim_X Y}$ in $U(1)$.
Thus we can equip it with the the canonical constant grading given by $-(\codim_Y X)\pi$ in $\R$.
\end{proof}

We will see later that with the canonical grading on $L_{Y,f}$, the Fukaya morphism complex
$\hom_{Fuk(T^*X)}(L_{X},L_{Y,f})$
has cohomology equal to the cohomology $H^*(Y)$ with its usual grading.
Here we have written $L_X$ for the zero section $T^*_X X$ with its canonical grading.


\subsubsection{Relative pin structure}
Recall first that the group $Pin^+(n)$ is the double cover of
$O(n)$ with center $\Z/2\Z \times \Z/2\Z.$\footnote{There is another
double cover $Pin^-(n)$ with center $\Z/4\Z.$} 
A pin structure on a Riemannian manifold $L$ is a lift
of the structure group of $TL$ to $Pin^+(n).$
The obstruction to a pin structure is the second
Stiefel-Whitney class $w_2(L)\in H^2(L,\Z/2\Z)$,
and choices of pin structures form a torsor over the group
$H^1(L,\Z/2\Z)$. 

A relative pin structure on a
submanifold $L \hookrightarrow M$ with background class $[w]\in H^2(M,\Z/2\Z)$
can be defined as follows.  Fix a {\v C}ech cocycle $w$ representing $[w]$,
and let $w|_L$ be its restriction to $L$. Then a pin structure on $L$ relative to $[w]$ can be
defined to be  
an $w|_L$-twisted pin structure on $TL$. Concretely, this can be represented by a $Pin^+(n)$-valued
\v Cech $1$-cochain on $L$ whose coboundary is $w|_L$. 
Such structures are canonically independent of the choice of \v Cech representatives.

\begin{rmk}
For a given background class $[w]$, 
choices of relative pin structures on $L$ form a torsor over the group
$H^1(L,\Z/2\Z)$. 
\end{rmk}

We check that our standard Lagrangians have canonical
relative pin structures with respect to a canonical universal background class. 
Recall that to a submanifold $Y\subset X$, and a defining function $m:X\to\R_{\geq 0}$
for the boundary $\del Y\subset X$, we have the standard Lagrangian 
$$
L_{Y,f}= T^*_Y X + \Gamma_{df} \subset T^*X
$$ 
where $f:Y\to \R$ is given by $f=\log m$.

\begin{prop}
The second Stiefel-Whitney class $w_2(L_{Y,f})\in H^2(L_{Y,f},\Z/2\Z)$
is the restriction of $\pi^*(w_2(X)).$
In fact, there is 
a canonical relative pin structure on $L_{Y,f}$ with background
class  $\pi^*(w_2(X)).$

\end{prop}

\begin{proof}
Since there is a canonical homotopy class of isotopies between
$L_{Y,f}\hookrightarrow T^* X$ and 
$T^*_Y X\hookrightarrow T^*X$, it suffices to check the
assertion for the latter. 
The metric provides a canonical isomorphism
between the restriction $T T^*_Y X|_Y$ and the restriction $TX|_Y$.
By functoriality, we have
the desired relative pin structure.
\end{proof}

\subsubsection{Definition of brane structures}
Finally, we have the definition of a brane structure on a Lagrangian.

\begin{dfn}[\cite{Seidel}] A {\em brane structure} $b$ on
a Lagrangian submanifold $L\subset T^*X$ is a pair $b = (\widetilde{\alpha},P)$ where
 $\widetilde{\alpha}: L \rightarrow \R$ is
a lift of the squared phase map,
and $P$ is a
relative pin structure on $L.$
\end{dfn}

We have seen in the above discussion that our standard Lagrangians come equipped with
canonical brane structures. We refer to a standard Lagrangian equipped with its
canonical brane structure as a {\em standard brane}.


\subsection{Definition of Fukaya category}

In this section, we define the Fukaya $A_\infty$-category $Fuk(T^*X)$.
General foundations are taken largely from \cite{Seidel}, and 
we restrict the discussion here to issues arising from noncompact Lagrangians.
We assume $X$ is a compact, Riemannian, real analytic manifold, and equip 
$T^*X$ with its canonical exact symplectic structure.
Throughout what follows, we fix an analytic-geometric category $\CC$
and assume all subsets are $\CC$-subsets unless otherwise noted.


\subsubsection{Objects}

Fix once and for all the canonical trivialization $\eta^2$ of the bicanonical bundle $\kappa^2$,
and the background relative pin class $\pi^*(w_2(X))$. All brane structures will be with 
reference to these fixed structures. 
The following should be considered a preliminary definition
until we discuss perturbations. 

\begin{dfn}
Objects of $Fuk(T^*X)$ are quadruples $({L},\mathcal E,b, \Phi),$ 
where ${L} \subset{T}^*X$ is an exact Lagrangian submanifold
such that $\ol L \subset  \ol T^*X$ is a $\CC$-subset,
$\mathcal E \rightarrow L$ is a vector bundle with flat connection,
$b = (\widetilde{\alpha},P)$ is a brane structure on $L$,
and $\Phi$ is a collection of perturbations to be explained below.
\end{dfn}

When circumstances are clear, we often refer to an object of $Fuk(T^*X)$ by its corresponding
support Lagrangian.
Given a submanifold $Y\subset X$ equipped with a local system $\CL_Y$, 
we refer to a standard Lagrangian $L_{Y, f}$ equipped
with the flat bundle $\CE= \pi^*(\CL_Y\otimes or_X \otimes or_Y^{-1})$,
and its canonical brane structure $b$ as a {\em standard object}.

We have defined the objects from the point of
view of the compactified cotangent bundle $\ol T^*X$ in order to give
a cleaner definition of the Lagrangians of interest.  Requiring
$\overline{L}$ to be a $\CC$-subset
of $\overline{T}^*X$ excludes various
types of behavior near infinity $T^\infty X$. For example, with our definition,
we can not have infinitely many intersection
points (as might occur for a helix on $T^*S^1$).  
Although we rule this out from the beginning, certain theories of the Fukaya category on $T^*X$ 
allow such behavior. 
Note that we use $\overline{T}^*X$ as a topological compactification, but
not as a symplectic compactification.
From the point of view of constructing moduli spaces (see below),
our Lagrangians are noncompact.


\subsubsection{Morphisms}
\label{homsec}
To define the morphisms between two objects, we need to choose
Hamiltonian isotopies to move their underlying Lagrangians so that they do not intersect
at infinity, and have transverse intersections in finite space. 
As usual, the
intersections will depend on the choice of isotopies, but in a homotopically manageable way. 
First, we explain here a broad class of isotopies which provide a consistent topological form
for the intersections of our Lagrangians. 
In the next section, we assume the existence of a more restricted class of
isotopies which guarantee that we may use moduli spaces of holomorphic disks 
to define composition maps.

Recall that a Hamiltonian function $H:T^*X\to \R$ is said to be {controlled} if there is a compact set $K\subset T^*X$ such that outside of $K$ we have $H(x,\xi)=|\xi|$. 
The corresponding Hamiltonian isotopy 
$\varphi_{H,t}:T^*X\to T^*X$ equals the normalized geodesic flow $\gamma_t$
outside of $K$. 
By Lemma~\ref{deform}, for Lagrangians $L_0,L_1\subset T^*X$, we may choose controlled
Hamiltonian functions $H_0,H_1$ and a fringed set $R\subset \R^2$
such that for $(\delta_0,\delta_1)\in R$, there is $r>0$ such that
$$
\varphi_{H_0,\delta_0}(\ol L_0)\cap 
\varphi_{H_1,\delta_1}(\ol L_1) \subset T^*X_{|\xi|<r},
$$
and the intersection is transverse.
Suppose we consider objects of $Fuk(T^*X)$ to come equipped with such data,
and that the brane structures and bundles are transported via the perturbations.
Then we may make the following definition.

\begin{dfn} For objects $L_0,L_1$ of $Fuk(T^*X)$, 
the space of morphisms is defined to be
$$
{hom}_{Fuk(T^*X)}(L_0,L_1) = \bigoplus_{p\in \varphi_{H_0,\delta_0}(L_0)\cap 
\varphi_{H_1,\delta_1}(L_1) }
{\mathcal Hom}(\mathcal E_0 \vert_{p},\mathcal E_1\vert_{p})[-\deg(p)].
$$
The integer
$\deg(p)$ denotes the Maslov grading, or index, of the linear
Lagrangian subspaces at the
intersection -- see Sections 11e-11g of \cite{Seidel}.
\end{dfn}

The differential on the complex of morphisms will be defined in Section \ref{compmaps}
below along with all of the higher composition maps.


\subsubsection{Holomorphic disks}
\label{hol disks}

The composition maps of the Fukaya $A_\infty$-category are defined
by counting points (with orientations) in appropriate moduli spaces of holomorphic maps
with respect to a compatible almost complex structure.
To ensure that the moduli spaces are well-behaved, one must consider further
perturbations, as described by Seidel \cite{Seidel}.
One must choose Floer perturbation data consisting of
a time-dependent Hamiltonian function and almost complex
structure deformation. 
One must also choose perturbation data on the Riemann surfaces to be mapped.

At the same time, we must check that the moduli spaces
are compact. 
This is delicate due to the fact that our Lagrangians are not necessarily compact.
For the case of closed Riemann surfaces
and no Lagrangians, if we consider surfaces intersecting a fixed compact set 
and with area less than a uniform upper bound, then tameness of the target
in the sense of \cite{Sikorav, ALP} ensures compact
moduli spaces. It is a standard fact that 
 $T^*X$ with its canonical symplectic form and the Sasaki almost complex structure $J_{Sas}$
  associated
to a Riemannian metric on $X$ is tame. To deal with the boundary case with
standard Lagrangians as boundary,
it will be useful to consider the conical almost complex structure $J_{con}$ introduced
in Section~\ref{warp}.
Recall that near infinity, the conical metric $g_{con}$ presents $T^*X$
as a cone over the unit (co-)sphere bundle $S^*X$.
It is straightforward to calculate explicitly that such a metric
has sectional curvature bounded from above
and injectivity radius bounded away from zero.
Thus $T^*X$ equipped with the conical almost complex structure $J_{con}$ 
is tame.

For the
Lagrangian boundary case, one typically imposes additional conditions
on the Lagrangian submanifolds themselves.  
The fact that our Lagrangians are exact means that none of the
complications from the bubbling of spheres will be present, but 
compactness still must be ensured.
To achieve compact moduli spaces,
there are many possible strategies involving assumptions adapted to
different situations.
The situation for standard branes
is robust: one can successfully apply a diverse host of
techniques including convexity statements, dilation arguments, and energy
bounds.
Thus in the context of this paper, our choice of a specific
uniform definition of $Fuk(T^*X)$ is largely aesthetic.
As with our use of Lagrangians lying in some analytic-geometric category,
we have decided upon the approach described below since it is particularly easy to work with in applications to mirror symmetry and representation theory.
We will apply the sufficient tameness conditions for compact moduli spaces 
derived in \cite{Sikorav, ALP}.  They require that (1) there exists $\rho_L>0$ such
that for every $x\in L$, the set of
points $y\in L$ with $d_{}(x,y) \leq \rho_L$ is
contractible, and that (2)
there exists $C_L$ giving a two-point condition
$d_L(x,y) \leq C_L d_{}(x,y)$  whenever $x,y\in L$ with $d_{}(x,y)<\rho_L$.

As the following example shows, standard Lagrangians
do not necessarily satisfy the above conditions with respect to the Sasaki metric $g_{Sas}$.

\begin{ex}
Consider any smooth curve $C\subset \R^2$ and its conormal bundle $T^*_C\R^2
\subset T^*\R^2$. If $C$ has an inflection point,
then $T^*_C\R^2$ is not tame
with respect to the Sasaki metric $g_{Sas}$.
Note that having an inflection point
is a generic circumstance.
\end{ex}

The reason for the above phenomenon is that the Sasaki metric $g_{Sas}$ is very
asymmetric: vertical directions
along co-sphere bundles grow with the radius,
while horizontal directions remain a fixed length.
We have introduced the conical metric $g_{con}$ to remedy the situation.
It is straightforward to check that for any compact submanifold $Y\subset X$, 
the conormal bundle $T^*_Y X \subset T^*X$ is tame with respect to
the conical metric $g_{con}$ (near infinity, the conormal is nothing more than a smooth subcone
of the ambient cone). While this fact is crucial, it is not the end of the story.
As the following example shows, standard Lagrangians
associated to submanifolds with singular boundaries are not necessarily tame
with respect to $g_{con}$.

\begin{ex}
Fix a point $p\in \R^2$ and a smooth closed curve $c:[0,\infty)\to \R^2$
with $c(0) = p$. 
Consider the complement $U = \R^2\setminus c$, and a standard
Lagrangian $L_U\subset T^*\R^2$ associated to $U$. 
If $c$ has non-vanishing curvature 
in a neighborhood of $p$, then $L_U$ will not be tame with respect to 
the conical metric $g_{con}$.
 Note that non-vanishing curvature is a generic circumstance. 
\end{ex}

Now to ensure that we have compact moduli spaces, we will give ourselves some added
flexibility and {assume}
that our Lagrangains come equipped with perturbations $\Phi$ making them tame
near infinity with respect to the conical metric $g_{con}$. 
For our intended applications, this assumption is easily verified and poses
no further restrictions.
(Furthermore, if one is willing to work with immersed but not necessarily
embedded Lagrangians, one could expect that there is no obstruction
to finding such perturbations. Since the foundations
of the Fukaya category of immersed Lagrangians are still not available in the literature,
we will not pursue this direction and insist for now that the perturbations exist.)

To make this precise, for a brane $L$ to define an object of $Fuk(T^*X)$,
we require the existence of the following further perturbation data.
First, by a perturbation of $L$, we mean a one-parameter
family of branes $L_t\subset T^*X \times \R$ such that $L_0 = L$
and for $|\xi|> r>0$ sufficiently large, the product map $L_t \to \R \times (r,\infty)$ given by 
the parameter $t$ and the length $|\xi|$ is a proper submersion. (This guarantees
that the family of branes has constant topology near infinity.)
Now fix a defining function $n:\ol T^*X\to \R_{\geq 0}$ for the closure $\ol L\subset \ol T^*X$,
and for $\varepsilon> 0$, let $N_{\vareps}(\ol L)\subset \ol T^* X$ be the open neighborhood 
$\ol T^*X_{n< \vareps}$. 
Given a brane $L$, we require the existence of a smooth perturbation $L_t$ such that
for all $\vareps>0$, the brane $L_{\vareps}$ is tame
with respect to the conical metric $g_{con}$
and lies in $N_{\vareps}(L)$.
It is worth pointing out that we do not insist that the family $L_t$ is uniformly tame with respect to $t$.
(In fact, if the original brane $L$ is not tame, then of course it is impossible
to find a uniformly tame perturbation.)
Similarly, we do not insist that the family $L_t$ extends to infinity.
Neither circumstance causes any harm.

\begin{lem}
Standard Lagrangians admit perturbations to Lagrangians tame with respect
to the conical metric $g_{con}$.
\end{lem}

\begin{proof} We will provide a concrete perturbation for standard branes.
The interested reader will note that 
the underlying construction is quite general and could be applied to many branes.
We will
move our brane to a new brane $L\subset T^*X$ such that in a neighborhood
of infinity, its closure $\ol L\subset \ol T^*X$ is diffeomorphic to the product
$L^\infty \times (r, \infty]$  of the boundary $L^\infty\subset T^\infty X$ with an interval.
Moreover,
near infinity the new brane $L$ will be uniformly close to the cone over its boundary $L^\infty$.
Since cones over compact submanifolds are tame with respect to the conical metric $g_{con}$,
this will immediately imply that $L$ is also tame. In fact, though it is unnecessary,
we could go one step further
and move $L$
so that near infinity it is equal to the cone over its boundary $L^\infty$.

With these general remarks in mind, let us turn to the case at hand of standard branes. 
Fix a submanifold $Y\subset X$ and
a defining function $m:X\to\R_{\geq 0}$
for the boundary $\del Y\subset X$, and consider
the standard Lagrangian $L_{Y,f}= \Gamma_{d f}|Y + T^*_{Y} X\subset T^*X|_Y\subset T^*X$ 
where as usual $f=\log m$.
Our first step will be to move $L_{Y,f}$ to a standard Lagrangian associated to 
a submanifold with smooth boundary.
For small $\eta>0$, 
choose an increasing function $b_{\eta}:\R\to \R$ satisfying $b_\eta(s) = s$,
for $s\geq \eta$, and $b_\eta(s) = 0$, for $s\leq 0$.
Consider the composition $m_\eta = b \circ m$ and the
submanifold $Y_\eta = Y_{m_\eta > 0} \subset Y$.
Then for all sufficiently small $\eta>0$, the boundary $\del Y_\eta \subset X$ is a smooth
submanifold. Our first step is to perturb $L_{Y,f}$ to the standard Lagrangian $L_{Y_\eta ,f_\eta}
=\Gamma_{d f_\eta}|Y_\eta + T^*_{Y_\eta} X \subset T^*X|_{Y_\eta}\subset T^*X$
where we have set $f_\eta = \log m_\eta$.

Let us assess what we have accomplished so far.
If $Y$ is open, then we claim that $L_{Y_\eta ,f_\eta}$ is tame.
To see this, observe that $\del Y_\eta$ is a smooth hypersurface with normal coordinate $m$.
Thus one can check directly in local coordinates that the conormal bundle
$T^*_{\del Y_\eta} X$ uniformly approximates $L_{Y_\eta ,f_\eta}$ near infinity.
If $Y$ is not open, then we must add a second step to our perturbation.
In this case, observe that $\ol Y_\eta$ is a smooth submanifold with boundary $\del Y_\eta$
with normal coordinate $m$.
Thus one can check directly in local coordinates that in
a neighborhood
of infinity, the closure $\ol L_{Y_\eta ,f_\eta}\subset \ol T^*X$ is homeomorphic to the product
$L_{Y_\eta ,f_\eta}^\infty \times (r, \infty]$ of the boundary $L_{Y_\eta ,f_\eta}^\infty\subset T^\infty X$
with an interval. Moreover, the boundary $L^\infty_{Y_\eta ,f_\eta}$
is a piece-wise smooth Legendrian (if $Y$ is closed, then it is in fact smooth; if $Y$ is not closed, 
it is a union of two submanifolds glued along their boundaries). Therefore we can perturb $L^\infty_{Y_\eta ,f_\eta}$ 
to a smooth nearby Legendrian $\CL^\infty$, and hence perturb $L_{Y_\eta ,f_\eta}$ to a nearby
Lagrangian which near infinity is uniformly approximated by the cone over $\CL^\infty$.
\end{proof}

Before continuing, let us make a couple of remarks.
First, a comment about the application of the above perturbations in the next section:
given a brane $L$ admitting such a perturbation, it follows that its image
under a controlled Hamiltonian isotopy does so as well.
Namely, we can simply conjugate the perturbation by the isotopy
since normalized geodesic flow is an isometry of the conical metric.
It is in this form that we will use the discussion of this section.
Second, although it is not necessary
for the aims of this paper, it is worth pointing out that one could allow
the flexibility of any compatible almost complex structure
as long as some $C^0$-bounds were obtained.
What we have described
is a simple concrete framework to deal with standard branes and objects
which arise in our intended applications.


\subsubsection{Composition maps}
\label{compmaps}

Now we are ready to define the composition maps of the Fukaya $A_\infty$-category
$Fuk(T^*X)$. 

Let $L_0,\dots,L_d$ be a finite collection of objects of $Fuk(T^*X).$
By Lemma~\ref{deform}, we may choose controlled
Hamiltonian functions $H_i:T^*X\to \R$, for $i=0,\ldots, d$,  and a fringed set $R\subset \R^{d+1}$
such that for $(\delta_0,\ldots,\delta_d)\in R$, there is $r>0$ such that
$$
\varphi_{H_i,\delta_i}(\ol L_i)\cap 
\varphi_{H_j,\delta_j}(\ol L_j) \subset T^*X_{|\xi|<r},
\mbox{ for $i\not =j$},
$$
and the intersections are transverse.
Then as discussed in the previous section, by assumption,
we may apply further small perturbations so that the resulting
branes are tame with respect to the conical almost complex structure.

We consider objects of $Fuk(T^*X)$ to come equipped with such data. 
Then by the results of \cite{Sikorav}, the moduli spaces
of holomorphic maps from 
a disk with a fixed number of marked boundary points taken to intersection points,
and boundary arcs taken to individual branes is compact.
In fact, by  \cite{Sikorav}, we have an explicit diameter bound on the image
of any such holomorphic map.
From this, one can verify that our perturbation framework can be handled
by standard techniques: the fringed set at infinity adds only a contractible
space of data to the usual Floer theory perturbations with
compact support, and one checks that the remaining
noncompact perturbations lead to canonical isomorphisms.
(See~\cite{Nad} for a
discussion of setting up the relevant continuation maps.) 
In summary, the usual methods of organizing the perturbation data 
to obtain a well-defined $A_\infty$-category apply (we refer the reader
to the careful account in \cite{Seidel}, Chapters 9 and 12).
For example, for a compactly supported time-dependent Hamiltonian $H_t$,
moduli spaces with moving boundary conditions
provide isomorphisms between the identity functor and the family of functors $\varphi_{H_t,T}$.
(See~\cite{Nad} for a discussion of which noncompactly supported
Hamiltonian isotopies provide isomorphisms.)

With the preceding in hand, we define the $A_\infty$-composition maps of $Fuk(T^*X)$
as usual
by their structure
constants: they count the signed number of holomorphic maps from 
a disk with $d+1$ marked boundary points with the appropriate boundary 
conditions.

\begin{dfn}
For $L_0,\dots,L_d$ objects of $Fuk(T^*X),$
one defines 
$$m^d_{Fuk(T^*X)}:
{hom_{Fuk(T^*X)}}(L_0,L_1)\otimes\dots\otimes{hom_{Fuk(T^*X)}}(L_{d-1},L_d)\rightarrow
{hom_{Fuk(T^*X)}}(L_0,L_d)
$$ 
as follows.  Consider elements 
$p_i\in {hom_{Fuk(T^*X)}}(L_i,L_{i+1}),$
for $i=0,\ldots,d-1$,
and $p_d\in hom_{Fuk(T^*X)}(L_0,L_d)$.
Then the coefficient of $p_d$ in 
$m^d_{Fuk(T^*X)}(p_0,\dots,p_{d-1})$
is the signed sum over holomorphic maps from a disk
with $d+1$ counterclockwise cyclically ordered
marked points mapping to the $p_i$ and corresponding
boundary arcs mapping to $L_{i+1}.$  Each map contributes
according to the holonomy of its boundary, where adjacent
components $L_i$ and $L_{i+1}$ are glued with $p_i.$
\end{dfn}


By Lemma~\ref{zeroarea},
there is no bubbling of spheres in the Fukaya category of an exact
symplectic manifold.
As a result, the Novikov coefficient rings employed to
account for all possible areas of maps are unnecessary, and we content
ourselves with simply counting the maps, with no weighting by areas.
This simplifies the isomorphism of different perturbation data,
since we need not keep track of the changes in the areas of disks
as the intersection points move.


\section{Embedding of standard objects}\label{a-infinity}

In this section we will construct an embedding of the
Morse $A_\infty$-category $Mor(X)$ into the Fukaya
$A_\infty$-category $Fuk(T^*X)$.
The embedding relies on rather delicate and
detailed perturbations of 
a collection of standard Lagrangians. 
After the necessary preparations, we will
be able to understand the moduli spaces of
holomorphic polygons bounding our perturbed Lagrangians
in terms of Morse theory via the theorem of Fukaya and Oh \cite{FO}.


\subsection{Preliminaries} We recall here some of our conventions and notations
concerning the geometry of the cotangent bundle $\pi:T^*X\to X$. Throughout what follows,
we identify $X$ with the zero section in $T^*X$.
We fix a Riemannian metric on $X$, and 
write $d_X(x,y)$ for the distance between points $x,y\in X$.

Under the metric identification $T(T^*X) \simeq T(TX)$,
the canonical one-form $\theta$ on $T^*X$ corresponds to
the geodesic vector field $v_\theta$. On the complement of the zero section $T^*X \setminus X$,
we have the normalized geodesic vector field
$\hat v_\theta= v_\theta/|v_\theta|$. It is the Hamiltonian vector field for the 
length function $H:T^*X \setminus X\to \R$ given by $H(x,\xi)=|\xi|$.
We write $\gamma_t:T^*X\setminus X\to T^*X\setminus X$ for the normalized geodesic flow for time $t$ associated to $\hat v_\theta$.
By definition, if we identify a covector $(x,\xi)\in T^*X$ with a vector $(x,v)\in TX$,
then we have the identity
$$
\gamma_t(x,v)=\exp_{x,t}(\hat v)_* (v)
$$
where $\hat v = v/|v|,$
the map $\exp_{x,t}:T_{x} X\to X$ denotes the exponential flow from the
point $x$ for time $t$, and the asterisk subscript indicates the
derivative (push-forward).
Note that for
$t$ sufficiently small -- for example, less than half the injectivity radius
of $X$ --
we have $\dist_X(\pi(x,\xi),\pi(\gamma_t(x,\xi)) = t.$

Given a stratification $\CS=\{S_\alpha\}$ of $X$, we define the associated conical Lagrangian 
$\Lambda_\CS\subset T^*X$ to be 
the union of conormals
$$
\Lambda_\CS=\cup_\alpha T^*_{S_\alpha}X.
$$
Given a second stratification $\CS'$ refining $\CS$, note the corresponding inclusion
$\Lambda_\CS\subset \Lambda_{\CS'}$.

Some usual notation:
given a space $Y$, a function $g:Y\to \R$, and $r\in \R$, we write $Y_{g=r}$ for the subset
$\{y\in Y| g(y)=r\}$, and similarly for inequalities. 


\subsection{Variable dilation}\label{dilation}
 We will need to dilate a standard Lagrangian
so that it is as close as we like to its associated conical Lagrangian.
To achieve this in a controlled fashion, we must consider two regions:
(1) a neighborhood of infinity where the Lagrangian is already close to
its associated conical Lagrangian, and (2) a compact region where
dilation of the standard Lagrangian is a Hamiltonian isotopy.

Let $U\subset X$ be an open submanifold, with closure $\ol U\subset X$,
and boundary $\del U=\ol U\setminus U$. 
Fix a defining function $m:X\to \R_{\geq 0}$ for the closed subset $X\setminus U$, and
let $f:U\to\R$ be the function $f=\log m$.
Let $L\subset T^*X$ be 
the standard Lagrangian given by the differential of $f$.

\begin{lem}\label{lemcomparison}
For any $\ell>0$, there is $\eta >0$ such that $L_{m<\eta}\subset L_{|\xi|>\ell}$.
For any $\eta>0$, there is $\ell>0$ such that $L_{m>\eta}\subset L_{|\xi|<\ell}$
\end{lem}

\begin{proof}
Immediate from the definitions.
\end{proof}

Choose a stratification of $X$ which refines the boundary $\del U\subset X$,
and let $\Lambda\subset T^*X$ denote the associated conical Lagrangian.
Choose a defining function $n:\ol T^*X\to \R_{\geq 0}$ for the closure $\ol\Lambda$,
and for any $\vareps\geq 0$,
let $N_\vareps(\Lambda)\subset T^*X$ denote the open neighborhood $T^*X_{n<\varepsilon}$
of $\Lambda$.

\begin{lem}
For any $\vareps>0$, there is $\ell>0$ such that 
$$
L_{|\xi|\geq \ell} \subset N_\vareps(\Lambda).
$$
\end{lem}

\begin{proof}
Recall that for a subset $Y\subset T^*X$, we write $Y^\infty$ for the intersection of the closure 
$\ol Y\subset \ol T^*X$ with the divisor at infinity
$T^\infty X\subset \ol T^*X$.
The assertion follows immediately from the inclusion $L^\infty\subset \Lambda^\infty$.
\end{proof}

Fix positive numbers $a<b\in (0,\infty)$, and choose an increasing function $d_{a,b}:\R\to \R$ satisfying the following
$$
d_{a,b}(r) = \left\{
{\begin{array}{cl}
 r, & \mbox{ for $r\geq \log b$,} \\
 \log\sqrt{ab}, & \mbox{ for $r\leq \log a$}.
\end{array}}\right.
$$

In order to dilate the Lagrangian $L$, we consider the Hamiltonian flow 
$\varphi_{D_{a,b},t}:T^*X\to T^*X$
generated by the function $D_{a,b}:T^*X\to \R$ defined by 
$$
D_{a,b}(x,\xi)=
\left\{
{\begin{array}{cl}
 -d_{a,b}(\log m(x)), & \mbox{ for $(x,\xi)$ with $m(x)\not = 0$,} \\
 -\log\sqrt{ab}, & \mbox{ otherwise}.
\end{array}}\right.
$$
The motion of the Lagrangian $L$ under the flow $\varphi_{D_{a,b},t}$
is given by the variable dilation
$$
\varphi_{D_{a,b},t}(L)=(1-t d_{a,b}'(\log m(x)))\cdot L.
$$
In particular, the Lagrangian $\varphi_{D_{a,b},t}(L)$ continues
to be a graph over $U$, and coincides with $L$ over $m\leq a$,
and with $(1-t)\cdot L$ over $m\geq b$.

\begin{lem}\label{neartoconical}
For any $\vareps>0$, there is $b>0$ and  $\delta>0$ 
such that for all $a',b'\in (0,b),$ $a'<b',$ and $\delta'>\delta$, we have
$$
\varphi_{D_{a',b'},\delta'}(L) \subset N_\vareps(\Lambda).
$$
\end{lem}

\begin{proof}
By the previous lemmas, we may choose $b>0$ so that 
$$
L_{m\leq b}
\subset N_\vareps(\Lambda).
$$
Since $\Lambda$ is conical,  $\varphi_{D_{a,b},t}$ preserves $N_{\vareps}(\Lambda)$.
Thus we need only choose $\delta>0$ so that  
$$
\varphi_{D_{a,b},\delta}(L_{m \geq b})\subset N_\vareps(\Lambda).
$$ 
But $L_{m \geq b}$ is a compact set, and $\Lambda$ contains the zero section of $T^*X$.
\end{proof}


\subsection{Separation}\label{separation}

We discuss here how to perturb a
standard Lagrangian near infinity. Namely, we show that near infinity we may separate
it from a conical Lagrangian without disturbing its structure elsewhere.

Let $U\subset X$ be an open submanifold, with closure $\ol U\subset X$,
and boundary $\del U=\ol U\setminus U$. 
Fix a defining function $m:X\to \R_{\geq 0}$ for $X\setminus U$, and
let $f:U\to\R$ be the function $f=\log m$.
Let $L\subset T^*X$ be 
the standard Lagrangian given by the graph of $df$.

Let $\CS=\{S_\alpha\}$  be any stratification of $X$, and let $\Lambda_\CS\subset T^*X$
be the corresponding conical Lagrangian 
$$
\Lambda_\CS=\cup_\alpha T^*_{S_\alpha} X.
$$
Note that we do not assume that $\CS$ has any relation to $U$ or its boundary $\del U$.

We say that $x\in X$ is a {\em $\Lambda_\CS$--critical point} of $m$ if we have $dm(x)\in\Lambda_\CS$.
Note that for $x\in U$ this is the same as $df(x)\in \Lambda_\CS$ since $df=dm/m$ and 
$\Lambda_\CS$ is conical. We say that $r\in \R$ is a {\em $\Lambda_\CS$--critical value} of $m$ 
if there is
a $\Lambda_\CS$--critical point $x\in X$ such that $r=m(x)$.

\begin{lem}\label{critvals}
There is $\eta>0$ so that there are no $\Lambda_\CS$-critical values of $m:X\to \R_{\geq 0}$ in the interval $(0,\eta]$.
\end{lem}

\begin{proof}
The $\Lambda_\CS$--critical values of $m$ form a discrete subset of $\R$.
\end{proof}

The following strengthening of Lemma~\ref{deform} will simplify
our perturbations, as we can choose the parameter $\eta>0$ 
to be independent of sufficiently small $\delta>0$.

\begin{lem}
There exist $\eta>0$ and $\delta>0$ 
such that for all $\delta'\in (0,\delta]$,
the normalized geodesic flow satisfies
$$
\gamma_{\delta'}(\ol L_{m\leq \eta})\cap \ol \Lambda_\CS =\emptyset.
$$
\end{lem}

\begin{proof}
We prove the assertion by contradiction. So suppose it were false. 

First, recall that for all sufficiently small $\delta>0$, we have seen in
Lemma~\ref{deform} that
$$
\gamma_\delta(L^\infty)\cap \Lambda_\CS^\infty = \emptyset.
$$
Thus by the previous lemma,
if the assertion were false,
then by the curve selection lemma, there exists a $\delta>0$,
and a subanalytic curve 
$$
\ell(t) = (x(t),\xi(t)):[0,\delta)\to L
$$
such that $m(x(t))\to 0$ as $t\to 0$, and (after a possible reparametrization) we have
$$
\gamma_{t}(\ell(t))\in T^*_{S_\alpha} X, \mbox{ for all $t\in (0,\delta)$},
$$
for some fixed
stratum $S_\alpha$ (we may fix $\alpha$
since there are finitely many strata).

Let $\kappa(t)=(y(t),\zeta(t))\in \Lambda_\CS$ denote the image curve $\gamma_{t}(\ell(t))$.
Again by definition, if we identify 
the covector $(x(t),\xi(t))\in T^*X$ with a vector $(x(t),v(t))\in TX$, then we have
$$
\zeta(t)(w)=
\langle \exp_{x(t),t} \hat v(t), w \rangle = 0,
\mbox{ for all $t\in (0,\delta)$, and $w\in T_{y(t)}S_{\alpha}$}.
$$

Let $x'(t)$ denote the tangent vector to the curve $x(t)$. 
Since $m(x(t))\to 0^+$ as $t\to 0^+$,
we have the inequality
$$
\xi(t)(x'(t))=\langle v(t), x'(t) \rangle = \langle \nabla f(x(t)),x'(t)\rangle
=\frac{d}{dt} f(x(t)) = \frac{1}{m(x(t))}\frac{d}{dt} m(x(t))> 0,
$$
for $t$ sufficiently small.

On the other hand, observe that $\dist_X(x(t),y(t)) = t$, so that 
$$
\frac{d}{dt}\dist_X(x(t),y(t))= 1.
$$
But in general, consider $x,y\in X$ connected by a geodesic
with tangent vector $v_x$ at $x$ and $v_y$ at $y$. Then for any curves
$x(t),y(t)$ in $X$, with $x(t_0)=x, y(t_0)=y$, we have
$$
\frac{d}{dt}\dist_X(x(t),y(t))|_{t=t_0}= \langle y'(t_0), v_y \rangle - \langle x'(t_0), v_x\rangle.
$$
In the case at hand,
$v_x = \hat{v}(t)$ and $v_y$ is its image under $\exp_{x(t),t}$.
But we have seen that
$
\langle y'(t_0), v_q \rangle = 0,
\mbox{ since $y'(t_0)\in T^*_{S_\alpha}X$,}
$
and also that
$
\langle x'(t_0), v_x\rangle > 0.
$
Thus we have
$$
\frac{d}{dt}\dist_X(x(t),y(t))|_{t=t_0}<0,
$$
and we have arrived at a contradiction.
\end{proof}

Fix positive numbers $k<\ell\in (0,\infty)$, 
and choose an increasing function $g_{k,\ell}:\R\to \R$
satisfying the following
$$
g_{k,\ell}(r) = \left\{
{\begin{array}{cl}
 r, & \mbox{ for $r\geq \ell$,} \\
 (k+\ell)/2, & \mbox{ for $r\leq k$}.
\end{array}}\right.
$$
Consider the Hamiltonian flow 
$\varphi_{G_{k,\ell},t}:T^*X\to T^*X$
generated by the function $G_{k,\ell}:T^*X\to \R$ defined by 
$$
G_{k,\ell}(x,\xi)=
 g_{k,\ell}(|\xi|).
$$
The flow $\varphi_{G_{k,\ell},t}$ is related to the normalized geodesic flow $\gamma_t$
by the formula
$$
\varphi_{G_{k,\ell},t}=\gamma_{g_{k,\ell}'(|\xi|)t}
$$
(recall $|\xi|$ is constant under $\gamma_t$).
In particular, $\varphi_{G_{k,\ell},t}$ is the identity when $|\xi|\leq k$, and is equal to $\gamma_t$
when $|\xi|\geq \ell$.

We have the following reformulation of the previous lemma.

\begin{lem}\label{sepnearinfinity}
There is $k>0$ and $\delta>0$ 
such that for all $\ell'> k' > k$ and $\delta'\in (0,\delta]$,
we have
$$
\varphi_{G_{k',\ell'},\delta'}(\ol L)\cap \ol \Lambda_\CS = L_{|\xi|< k'}\cap \Lambda_\CS
$$
\end{lem}

\begin{proof}
Immediate from the previous lemma and Lemma~\ref{lemcomparison}.
\end{proof}

\subsection{Perturbations}\label{perturbation section}

We are now ready to describe how to perturb a collection of standard Lagrangians.
In what follows, let $i$ denote an element of the index set ${\Z/(d+1)\Z}$.

Let $U_i\subset X$ be an open submanifold, with closure $\ol U_i\subset X$,
and boundary $\del U_i=\ol U_i\setminus U_i$. 
Fix a defining function $m_i:X\to \R$ for the closed subset $X\setminus U_i$, and
let $f_i:U_i\to\R$ be the function $f_i=\log m_i$. 
Let $L_i\subset T^*X$ be
the standard Lagrangian 
given by the graph of the differential of $f_i$.
Choose a stratification of $X$ which refines the boundary $\del U_i\subset X$,
and let $\Lambda_i\subset T^*X$ denote the associated conical Lagrangian.

We will apply a sequence of Hamiltonian perturbations
to the Lagrangians $L_i$ to put them in a good position. In order to satisfy
the definition of the Fukaya category, the perturbations
must be positive normalized geodesic flow near infinity. Furthermore, the amounts
$\delta_i>0$ of 
normalized geodesic
flow with which
we move the $L_i$ near infinity must satisfy
$$
(\delta_0,\dots,\delta_d) \in R \subset \R^{d+1},
$$
where $R$ is some fringed set.
Because of this requirement,
we will work backwards through the collection perturbing the Lagrangians in the order
$L_d,\ldots, L_0$.

At the $i$th stage,
each of our perturbations will consist of two steps.
(1) We will first variably dilate $L_i$ so
that it becomes 
arbitrarily close
to its associated conical Lagrangian $\Lambda_i$. 
(2) We will then gently perturb it near infinity in the direction
of positive geodesic flow.
The first step will have three effects: (a) all of the intersections of the resulting Lagrangian
with the previously perturbed Lagrangians 
will be near the zero section;  (b) the height of the resulting Lagrangian will be
less than that of the previously perturbed Lagrangians along certain critical contours;
(c) intersections of the resulting Lagrangian
with the associated conical Lagrangian of yet-to-be perturbed Lagrangians 
will be near the zero section.
The second step will ensure that the first step is effective.

For each Lagrangian $L_i$, we organize the discussion of its perturbation into four parts:
\begin{itemize}
\item
(Intersections) We first collect the other Lagrangians whose intersections with $L_i$ must
be either dilated close to the zero section, or perturbed away near infinity.
\item
(Dilation) The variable dilation to the zero section. 
\item (Separation) The small perturbation near infinity.
\item
(Conclusion) We finally organize the result so that we may proceed to the next Lagrangian.
\end{itemize}

Throughout, we fix a positive number $h>0$.
We begin with the last Lagrangian in the collection.\\

$L_d:$ (Intersections)
Since there are no previously perturbed Lagrangians,
our aim here is simpler than in general. 
Let $\Lambda_{\leq d}\subset T^*X$ be the conical Lagrangian
$$
\Lambda_{\leq d}=\cup_{j\leq d} \Lambda_j.
$$
To guarantee that intersections
with the yet-to-be perturbed Lagrangians $L_j$, for $j<d$, can be dilated close to the zero section $T_X^*X$, we
must dilate the intersection $L_d\cap \Lambda_{\leq d}$ close to $T_X^*X.$

(Dilation) 
By Lemma~\ref{critvals}, the intersection
$L_d\cap \Lambda_{\leq d}$ is compact. 

Therefore we may choose $\sigma_d>0$
such that the standard dilation satisfies
$$
(\sigma_d\cdot L_d)\cap  \Lambda_{\leq d}\subset N_{h/2}(T_X^*X),
$$
By compactness of the intersection and Lemma~\ref{lemcomparison}, 
we may choose $\eta_d>0$ so that 
$$
(\sigma_d\cdot L_d)\cap \Lambda_{\leq d}\subset T^* X_{m_d > \eta_d}.
$$
To truncate this dilation near infinity, choose
positive numbers $a_d<b_d\in (0,\eta_d)$
and a Hamiltonian function $D_{a_d,b_d}$ as in Lemma~\ref{neartoconical}.
The resulting variably dilated Lagrangian satisfies
$$
\varphi_{D_{a_d,b_d},\delta_d}(L_d)\cap \Lambda_{\leq d}\subset T^*X_{m_d>\eta_d,|\xi|<h/2}.
$$

(Separation) Next, we apply Lemma~~\ref{sepnearinfinity} to the dilated Lagrangian 
$\varphi_{D_{a_d,b_d},\delta_d}(L_d)$ and the conical Lagrangian $\Lambda_{\leq d}$. 
Let $M_d$ be the maximum of the length $|\sigma_d\cdot L_d|$
in 
the region $m_d\geq\eta_d$. 
We may choose 
$
k_d> \max\{h, M_d\}
$ 
and a Hamiltonian function $G_{k_d,\ell_d}$ such that for some $\epsilon_d>0$ 
the corresponding perturbation satisfies
$$
\varphi_{G_{k_d,\ell_d},\delta_d}(\varphi_{D_{a_d,b_d},\delta_d}(\ol L_d))\cap \ol \Lambda_{\leq d} 
=(\varphi_{D_{a_d,b_d},\delta_d}(\ol L_d))_{|\xi|< k_d}\cap \Lambda_{\leq d}.
$$

(Conclusion) 
We set
$$
\tilde L_d =\varphi_{G_{k_d,\ell_d},\delta_d}(\varphi_{D_{a_d,b_d},\delta_d}( L_d))\subset T^*X,
$$
$$
\tilde U_d = X_{m_d>\eta_d}\subset X,
$$
$$
\tilde \Gamma_d = (\tilde L_d)_{m_d>\eta_d, |\xi| < k_d}\subset \tilde L_d.
$$
Note that $\tilde \Gamma_d$ is a graph over $\tilde U_d.$
By construction, we have
$$
\tilde \Gamma_d= \sigma_d\cdot (L_d)_{m_d>\eta_d},
$$
$$
\ol{\tilde L_d}\cap  \ol \Lambda_{\leq d} ={\tilde \Gamma_d}\cap  \Lambda_{\leq d} \subset N_{h/2}(T_X^*X).
$$
\medskip

At an arbitrary step, we proceed as follows.
\\

$L_i$: (Intersections)
Let $\tilde L_{>i}\subset T^*X$ be the union of the previously perturbed Lagrangians
$$
\tilde L_{>i}=\cup_{j>i} \tilde L_j.
$$
We would like to dilate the intersection
with $\tilde L_{>i}$ close to the zero section $T_X^*X$.
We will not be able to move this intersection closer than the intersection $\tilde L_{>i} \cap \Lambda_i$, but
at least this intersection is already close by induction.

Let $\Lambda_{\leq i}\subset T^*X$ be the conical Lagrangian
$$
\Lambda_{\leq i}=\cup_{j\leq i} \Lambda_j.
$$
To guarantee that intersections
with the yet to be perturbed Lagrangians $L_j$, for $j<i$, can be dilated close to
$T_X^*X$, we
must dilate the intersection $L_i\cap \Lambda_{\leq i}$ close to $T_X^*X.$

Let $\Lambda_{> i}\subset T^*X$ be the union of conormals
$$
\Lambda_{> i}=\cup_{j> i} T^*_{\del \tilde U_j} X.
$$
To guarantee that there is not unmanageable behavior along the boundaries
of the previously defined open sets $\tilde U_j$, for $j>i$, 
we must dilate the intersection $L_i\cap \Lambda_{> i}$ close to $T_X^*X$.

We set
$$
\Lambda_{[i]}= \Lambda_{\leq i}\cup \Lambda_{> i}.
$$

(Dilation) 
By induction, we have
$$
\ol {\tilde L}_{>i} \cap \ol\Lambda_i
\subset N_{h/2}(T_X^*X).
$$
Since $L_i^\infty\subset \Lambda_i^\infty$, it follows that
$$
L^\infty_i \cap {\tilde L^\infty}_{>i} =\emptyset,
$$
and so the intersection $L_i\cap \tilde L_{>i}$ is also bounded.
In addition, by Lemma~\ref{critvals}, the intersection
$
L_i\cap \Lambda_{[i]}
$
is compact.

Therefore we may choose $\sigma_i>0$ such that the standard dilation satisfies
$$
(\sigma_i\cdot L_i)\cap \tilde L_{>i}\subset N_{h}(T_X^*X),
$$
$$
(\sigma_i\cdot L_i)\cap  \Lambda_{[i]}\subset N_{h/2}(T_X^*X).
$$
Furthermore, for $\sigma_i>0$ sufficently small,
we may arrange for 
$$
|\sigma_i\cdot L_i| < \min_{j>i}{M_j},
\mbox{ above the compact set 
$(\sigma_i\cdot L_i)\cap \Lambda_{>i}$,}
$$
where $M_j$ denotes 
the maximum of the length $|\sigma_i\cdot L_i|$ 
in the region $m_i\geq\eta_i$.
By compactness of the intersections and Lemma~\ref{lemcomparison}, 
we may choose $\eta_i>0$ so that 
$$
(\sigma_i\cdot L_i)\cap (\tilde L_{>i}\cup \Lambda_{[i]})\subset T^* X_{m_i > \eta_i}.
$$
To 
truncate this dilation near infinity, choose
positive numbers $a_i<b_i\in (0,\eta_i)$
and a Hamiltonian function $D_{a_i,b_i}$ as in Lemma~\ref{neartoconical}.
The resulting variably dilated Lagrangian
$
\varphi_{D_{a_i,b_i},\delta_i}(L_i)
$
satisfies all of the properties derived above for 
the standard dilated Lagrangian $\sigma_i\cdot L_i$.

(Separation) Recall that we have 
$$
L^\infty_i \cap {\tilde L^\infty}_{>i} =\emptyset.
$$
We apply Lemma~\ref{sepnearinfinity}
to the dilated Lagrangian 
$\varphi_{D_{a_i,b_i},\delta_i}(L_i)$ and the conical Lagrangian $\Lambda_{[i]}$.
Let $M_i$ be as above.
We may choose 
$
k_i> \max\{h, M_i\}
$ and a function $G_{k_i,\ell_i}$ and an 
$\epsilon_i>0$ 
such that 
$$
\varphi_{G_{k_i,\ell_i},\delta_i}(\varphi_{D_{a_i,b_i},\delta_i}(\ol L_i))\cap \ol \Lambda_{[i]} 
= 
(\varphi_{D_{a_i,b_i},\delta_i}(L_i))_{|\xi|< k_i}\cap \Lambda_{[i]}.
$$

(Conclusion) 
We set
$$
\tilde L_i =\varphi_{G_{k_i,\ell_i},\delta_i}(\varphi_{D_{a_i,b_i},\delta_i}( L_i))\subset T^*X,
$$
$$
\tilde U_i = X_{m_i>\eta_i}\subset X,
$$
$$
\tilde \Gamma_i = (\tilde L_i)_{m_i>\eta_i, |\xi| < k_i}\subset \tilde L_i.
$$
Note that $\tilde \Gamma_i$ is a graph over $\tilde U_i.$
By construction
$$
\tilde \Gamma_i= \sigma_i\cdot (L_i)_{m_i>\eta_i}  
$$
$$
\ol{\tilde L_i}\cap  \ol {\tilde L_{j}}  
={\tilde \Gamma_i}\cap   {\tilde \Gamma_{j}}  
\subset N_{h}(T_X^*X), \mbox{ for all $j>i$},
$$
$$
\ol{\tilde L_i}\cap  \ol \Lambda_{[i]} ={\tilde \Gamma_i}\cap  \Lambda_{[i]} 
\subset N_{h/2}(T_X^*X),
$$
$$
|\tilde \Gamma_i| < |\tilde \Gamma_j|
\mbox{ wherever
$\tilde \Gamma_i\cap T^*_{\del \tilde U_j} X$, for all $j>i$}.
$$
\newline

By following this procedure, we arrive at the following.

\begin{prop}
\label{pertout}
The collection of Lagrangians $\tilde L_i$, graphs $\tilde \Gamma_i$,
and open sets $\tilde U_i$ satisfies the following.
\begin{enumerate} 
\item $\ol{\tilde L}_i\cap \ol{\tilde L}_j = \tilde \Gamma_i \cap \tilde \Gamma_j$, for $i\not = j$.
\item $|\xi|^2$ has no critical points on $\tilde L_i\setminus \tilde \Gamma_i$.
\item $(\tilde U_i, \tilde \Gamma_i)$ form a transverse collection in $Mor(X)$ --
see Definition~\ref{transversedef}.
\end{enumerate}

\end{prop}

\begin{proof}
The last assertion is the only part left to check. By construction,
the collection of boundaries $\del\tilde U_i$ are transverse.
To see $(\tilde U_i, \tilde \Gamma_i)$ is transverse, we need only check
that there is a metric for which the corresponding difference vector fields
point in the appropriate inward and outward directions. Such a metric may be
constructed locally wherever the level sets of the defining functions are transverse.
By construction, at any places where transversality fails, the relative sizes of the vector
fields have been arranged to allow for a metric to be constructed.
\end{proof}


\subsection{Relation to Morse theory}\label{relation to morse theory section}

The PSS isomorphism refers to an equivalence between Floer
homology and singular homology, and appears in both Hamiltonian
and Lagrangian Floer theory -- see Section 3 of \cite{albers} for a
recent discussion.
In the context of Lagrangian graphs in the cotangent bundle of a compact manifold,
Fukaya and Oh \cite{FO} extended this to an identification of 
the Morse and Fukaya $A_\infty$-categories by establishing
an oriented diffeomorphism of the moduli
spaces of gradient trees and holomorphic polygons involved in the definition 
of the higher composition maps.
In the local setting of graphs over open sets with transverse boundaries,
Kasturirangan and Oh~\cite{KO1, KO2} prove an equivalence of the Morse and Floer chain complexes.
In this section,
we adapt the approach of Fukaya and Oh
to prove an $A_\infty$-equivalence
of Morse and Fukaya $A_\infty$-categories which include all standard objects.
To do this, we first recall
the theorem of Fukaya and Oh in its original
form (with notation modified to agree with ours),
then adapt our situation to be able to apply their constructions.

\begin{fothm}[\cite{FO}]  Let $(X,g)$ be a Riemannian
manifold, and let $J_{Sas}$ be the canonical (Sasaki) almost complex structure on
$T^*X.$  Let ${f} = (f_0,...,f_d)$ be a generic collection
of functions on $X$, and let 
${\Gamma} = (\Gamma_{df_0},...,
\Gamma_{df_d})$ be the graphs of their differentials. 
For sufficiently
small $\epsilon>0,$ 
there is an oriented diffeomorphism
between the Morse moduli space of gradient trees of ${f}$
and the Fukaya moduli space of pseudoholomorphic
disks (with respect to $J_{Sas}$)
bounding the Lagrangians $\epsilon{\Gamma}.$
\end{fothm}

Recall that Proposition~\ref{pertout} of the preceding section
provides, starting from a collection $L = (L_0,\ldots,L_d)$ of standard objects of $Fuk(T^*X)$,
a perturbed collection
$
{\tilde L}=(\tilde L_{0},\ldots, \tilde L_{d}).
$
Above the open set $\tilde U_i = X_{m_i>\eta_i}$, the perturbed object $\tilde L_i$ results from
dilating the original object $L_i$, and thus in particular remains a graph
$$
\tilde \Gamma_i = (\tilde L_i )_{m_i>\eta_i}= \sigma_i \cdot(L_i )_{m_i>\eta_i}.
$$
Furthermore, all intersection points of the compactifications $\ol {\tilde L_i}$ 
occur among the $\tilde \Gamma_i$,
and $|\xi|^2$ has no critical points on the complements $\tilde L_i \setminus \tilde \Gamma_i$.
Finally, in the category $Mor(X)$, we have a transverse collection of objects
$$
{\tilde \fU}=((\tilde U_0,\tilde f_0),\ldots,(\tilde U_d,\tilde f_d))
$$
where the graph of the differential $d \tilde f_i$ is precisely
$\tilde \Gamma_i$. 
In what follows, we write $\tilde \Gamma$ for the collection of partial graphs
$(\tilde\Gamma_0, \ldots, \Gamma_d)$.

\medskip

We can not apply the Fukaya-Oh theorem directly to the above situation for several reasons.
First, the perturbed Lagrangians $\tilde L_i$ are noncompact and no longer graphs.
Moreover, as described in the previous section, we need to consider the conical
almost complex structure $J_{con}$ which is only equal to the Sasaki almost complex structure 
$J_{Sas}$
near the zero section. 
Finally, the functions $\tilde f_i$ and corresponding graphs $\tilde \Gamma_i$ 
are defined only over the open sets $\tilde U_i$.
Instead, 
we pursue the following strategy. Let us restrict
our attention to the collection of bounded but partial graphs $\tilde \Gamma$. 
For small enough $\epsilon >0$,
the local nature of the Fukaya-Oh theorem will give an identification of Fukaya and Morse
moduli spaces
for the dilated collection $\epsilon \tilde\Gamma$. Here we are using
that the conical almost complex structure $J_{con}$ is equal to the Sasaki
almost complex structure $J_{Sas}$ near to the zero section. With this understood, 
we need only show that
the Fukaya moduli spaces 
for $\epsilon \tilde\Gamma$ continue to calculate an $A_\infty$-structure quasi-isomorphic 
to that of the original unbounded but complete collection $\tilde L$.
This comes down to showing compactness of the moduli spaces for $\epsilon \tilde\Gamma$
(and similarly, compactness of the
moduli spaces providing continuation maps) as we vary the dilation parameter~$\epsilon$.

Details of this approach follow below. Throughout we work with the conical 
almost complex structure $J_{con}$ and corresponding conical metric $g_{con}$
associated to a Riemannian metric on $X$. Our constructions will take place near the
zero section where these structures agree with the respective Sasaki structures.

\medskip

(1) (Area bounds) 
First, 
choose a small $\eta'_i>\eta_i$, and
consider the level-set 
$$
X_{m_i = \eta_i'}\subset X.
$$
Choose a small $\delta_{i,h}>0$, and define the annulus-like open set
$$
S_i \subset X
$$
to consist of all points whose distance to the level-set $X_{m_i = \eta_i'}$ is less than $\delta_{i,h}$.
Consider the annulus-like partial graph
$$G_i = \tilde \Gamma_i \cap \pi^{-1}(S_i).$$
Finally, choose a very small $\delta_{i, v} >0$, and define the tube-like open set
 $$
T_i \subset T^* X
$$
to be the union of 
the vertical balls
$B^v_{\delta_{i,v}}\subset T^*X$ of radius $\delta_{i, v}$ centered at points of $G_i$.
Here by the vertical ball $B^v_{\delta_{i,v}}$ around a covector $(x,\xi)\in T^*X$, 
we mean the ball in the fiber $T^*_x X$ centered at $\xi$.

By construction, for small enough $\delta_{i,h}>0$,
the boundary of the partial graphs
$
G_i
$
decomposes as a disjoint union of manfolds
$$
\partial G_i = H^{in}_i \cup H^{out}_i
$$
where $m_i|_{H^{in}_i }> \eta_i'$ and $m_i|_{H^{out}_i} < \eta_i'$.
Furthermore, for small enough $\delta_{i,v}>0$,
there will be no interaction among the tubes:
$$
T_i \cap T_j = \emptyset,
\mbox{ for $i\not = j$}.
$$

Our aim is to use the relatively compact region $T_i$ to construct an area bound
on holomorphic disks.
An important wrinkle is that we would like the bound to behave well with respect
to dilations of $\tilde L_i$ together with $T_i$ towards the zero section. 
The precise statement we need is contained in the following
monotonicity bound from \cite{Sikorav}.

\begin{lem}\label{second area bound}
There exist
constants $R_i, a_i > 0$
such that the following holds for any $0< \epsilon\leq 1$.

Consider a holomorphic map
$u: (D,\partial D)\rightarrow (\epsilon T_i, \epsilon G_i)$.
Then for any ball $B_{r}\subset \epsilon T_i$
of radius $ r<\epsilon R_i$ such that $u(D)$ contains the center of $B_{ r}$, we have
$$
{\rm Area}(u(D) \cap B_{ r}) > a_i  r^2.
$$
\end{lem}

\begin{proof}
For fixed $0<\epsilon\leq 1$, the assertion
follows from Proposition 4.7.2 of \cite{Sikorav}
(note that Sikorav's proofs of 4.3.1(ii) and 4.7.2(ii) are entirely local --
only the bounding constants are global in nature).

What remains is the assertion that the bound can be achieved
uniformly with respect to dilation by $\epsilon$.
But the family of graphs $\epsilon G_i$ extends
to a compact family including $\epsilon = 0$ where we simply take the zero section itself.
Thus all of the controls on the geometric complexity required to apply the 
area bound of~\cite{Sikorav} can be achieved uniformly
with respect to $\epsilon$.
\end{proof}

The above bounds for $\epsilon =1$ allow us to 
control where disks with boundary along $\tilde L$ can go.
More precisely, 
we can fix a small
radius $0<\mathfrak r_i<R_i$ such that any ball $B_{\mathfrak r_i}$ of radius $\mathfrak r_i$ 
centered at a point of $(G_i)_{m_i=\eta'}$ fits inside $T_i$.
Then by variably dilating each $\tilde L_i$
while fixing $T_i$, we can
 arrange that any 
disk along the resulting collection 
with a fixed number $d$ of marked points
has area less than the minimum of the above bounds 
$a_i \mathfrak r_i^2$. 
Thus the boundaries of the disks can not pass through the regions $G_i$
and so must lie on the variably dilated $\tilde \Gamma_i$.

\medskip

(2) (Uniform dilation) Next, we take advantage of the homogeneity of the above area bounds
to see that we need not restrict ourselves to variable dilations of the $\tilde L_i$, but can in fact
dilate the entire $\tilde L_i$ uniformly. The main point is that the argument
of the preceeding paragraph, which applies to $\epsilon = 1,$
is robust enough
to allow the $T_i$ to be dilated to $\epsilon T_i$ as well.

To begin, we refine some of the choices made in the previous step. First, fix
any $0<\delta'_{i,h}< \delta_{i,h}$, and consider the smaller annulus-like region
$$
S'_i \subset X
$$
consisting of all points whose distance to the level-set $X_{m_i = \eta_i'}$ is less than $\delta'_{i,h}$.
Similarly, consider the annulus-like partial graph
$$G'_i = \tilde \Gamma_i \cap \pi^{-1}(S'_i).$$

For any such $\delta'_{i,h}$, we can find $0<\mathfrak r_i< R_i$ such that for any ball $B_{\mathfrak r_i}
\subset T^*X$
of radius $\mathfrak r_i$ centered at any point of $G'_i$, we have
$$
B_{\mathfrak r_i} \subset T_i.
$$
Fix such a radius $\mathfrak r_i$ and
also fix the area bound $a_i$ given by Lemma \ref{second area bound}.
As explained above, after variably dilating each $\tilde L_i$ while fixing the $T_i$, we may assume that
the area of any disk along the resulting collection satisfies
$$
{\rm Area}(u(D)) < a_i  \mathfrak r_i^2.
$$

Now, consider the absolute dilation $\epsilon \tilde L$ for any $0< \epsilon\leq 1$.
We claim that for any holomorphic map 
$u: (D,\partial D)\rightarrow (T^*X, \epsilon \tilde L)$ 
with fixed number of marked points $d$,
the boundary $u(\partial D)$
must in fact lie in $\epsilon \tilde\Gamma$. 
Suppose otherwise, and consider a boundary path $u(C)$ in $\epsilon G_i$
traversing from $\epsilon H_i^{in}$ to $\epsilon H_i^{out}$.
By construction, through subdividing $S'_i$ and considering
the induced subdivision of $\epsilon T_i$, we can find
$1/\epsilon$ disjoint balls $B_{\epsilon \mathfrak r_i}\subset \epsilon T_i$
of radius $ \epsilon \mathfrak r_i$ such that the path contains the centers of
the $B_{\epsilon \mathfrak r_i}$.
Thus the lemma gives the area bound
$$
{\rm Area}(u(D) \cap \epsilon T_i) > \epsilon a_i  \mathfrak r_i^2.
$$
Since the possible area of disks with boundary on the collection also scales linearly with
the dilation paramter $\epsilon$
({\it cf.} Lemma \ref{zeroarea}), 
we conclude that the initial area bound
for $\epsilon =1$ implies that the boundary $u(\partial D)$
must lie in $\epsilon \tilde\Gamma$.

\medskip

(3)  (Application of Fukaya-Oh theorem)
Finally, we apply the Fukaya-Oh theorem. The key observation to make
is that the proof of Fukaya-Oh given in~\cite{FO} is local
in the following sense. Given a gradient tree in $X$ for the functions
$f$,
the corresponding holomorphic disk in $T^*X$ for the graphs $\epsilon \Gamma$ will be 
in a neighborhood of the gradient tree such that the size of the neighborhood
goes to zero as $\epsilon$ goes to zero. And conversely, any such holomorphic disk
will arise in this way.

More precisely, starting from a gradient tree,
Fukaya and Oh first construct an approximate holomorphic disk
$w_\epsilon$.  The distance between
the gradient tree and $w_\epsilon$ goes to zero as $\epsilon \rightarrow 0.$
Next, an actual solution is proven to exist nearby in the $L^\infty$-topology. 
The actual solution (in the notation of \cite{FO}) has the
form $\exp_{w_\epsilon}(Q\eta)$.
The point is that by Theorem 9.1 of~\cite{FO}, for all $\delta>0$
there exists
an $\epsilon > 0$ such that for all $0 < \epsilon' < \epsilon,$
one has $\Vert Q\eta\Vert_{L^\infty} < \delta.$
In other words, the distance between the actual solution and
the gradient tree can be made arbitrarily small as long as $\epsilon$ goes to zero.

We conclude that by dilating our partial graphs
$\tilde\Gamma_i$ over the open sets
$\tilde U_i$ uniformly close to the zero section,
any gradient tree for the collection $\tilde f$
will correspond to a holomorphic disk with boundary
on the $\tilde \Gamma_i$. 
Finally, the area bounds
from the previous step implies such holomorphic disks
are the only ones to be considered.
We thus have arrived at our desired result.

\begin{thm}
There is an $A_\infty$-quasi-equivalence between $Mor(X)$ and
the full subcategory of
$Fuk(T^*X)$ generated by the standard objects $L = \Gamma_{df}$
over open sets $U\subset X$,
where $f:U\to \R$ is given by $f=\log m$, and $m:X\to \R_{\geq 0}$
is a defining function for the complement $X\setminus U$.
\end{thm}


\section{Arbitrary standard objects}\label{arbitrary}

For future applications, it is useful to know where the embedding takes other objects 
and morphisms. In particular, we would like to know not only where it takes standard sheaves
on open submanifolds, but also standard sheaves on arbitrary submanifolds. 
As discussed in the introduction, one approach
to this problem is to express standard sheaves on arbitrary submanifolds in terms of standard sheaves
on open submanifolds, and then to check what the relevant distinguished
triangles of constructible sheaves look like under the embedding. This requires
identifying certain cones in the Fukaya category with symplectic surgeries. 
Rather than taking this 
route,
we will instead show in this section that we may explicitly extend the domain of the embedding to include
standard sheaves on arbitrary submanifolds and morphisms between them.

We will follow very closely the steps used to define the embedding
in the preceding sections.
First, we will interpret the dg category $Sh(X)$ of constructible sheaves in terms of a category
$Sub(X)$ whose objects are submanifolds (equipped with certain defining functions)
and whose morphisms are complexes of relative de Rham forms (on certain open submanifolds with hypercorners).
Next, we will interpret the category $Sub(X)$ in terms of an extended version
of the category $Mor(X)$ built out of Morse theory. Finally, we will explain how the work of Fukaya-Oh
may be adapted to identify $Mor(X)$ with a full subcategory of the Fukaya category $Fuk(T^*X)$.
Because of the amount of overlap with the preceding sections, we will only explain 
the new wrinkles which arise and not repeat all details.

Before continuing,
we state here where the embedding takes the standard sheaf $i_*\CL_Y$
associated to a local system $\CL_Y$ on an arbitrary submanifold $i:Y\hookrightarrow X$.
Suppose that we are given a defining function $m:X\to \R_{\geq 0}$ for the boundary $\del Y\subset X$.
Recall that we define the standard Lagrangian $L_{Y,m}\subset T^*X$ to be the fiberwise sum
$$
L_{Y,m} =T^*_Y X + \Gamma_{d\log m}
$$
where $T^*_YX\subset T^*X$ is the conormal bundle to $Y$, and $\Gamma_{d\log m}\subset T^*X$
is the graph of the differential of $\log m$.
As explained in Section~\ref{branes},
it comes equipped with a canonical brane structure $b$, along with a 
flat bundle $\CE=\pi^*(\CL_Y\otimes or_X \otimes or_Y^{-1})$.
We write $L_{Y,m,\CL_Y}$ for the corresponding object of $Fuk(T^*X)$.
The main consequence of this section is the following.

\begin{thm}
Under the quasi-embedding $Sh(X)\hookrightarrow \tw Fuk(T^*X)$, the image of
the standard sheaf $i_*\CL_Y$ is canonically isomorphic to the standard brane $L_{Y,m,\CL_Y}$.
\end{thm}

In what follows, we limit the discussion
to the case of trivial local systems since the arbitrary case is no more difficult.
This will help streamline the exposition -- for example,
we write $L_{Y,m}$ for the object $L_{Y,m,\CL_Y}$ when $\CL_Y$ is trivial.


\subsection{Submanifold category} In Sections~\ref{open} and~\ref{smooth}, we introduced
the dg category $Open(X)$.
The results of this section generalize that discussion.

We define a dg category $Sub(X)$ as follows. The objects of $Sub(X)$ are triples
$(Y,m,n)$ where $Y\subset X$ is a submanifold, $m:X\to \R_{\geq 0}$ is a defining function
for its boundary $\del Y\subset X$, and $n:X\to \R_{\geq 0}$
is a defining function
for its closure $\ol Y\subset X$. To define
the complex of morphisms from an object $\fY_0=(Y_0,m_0,n_0)$
to an object $\fY_1=(Y_1,m_1,n_1)$, we introduce some perturbations.
It will be clear that the choices range over a contractible set, and that they can be made
compatibly for any finite collection of objects.
We will use the following general statement repeatedly.

\begin{lem}\label{tubelemma}
Let $\fY=(Y,m,n)$ be an object of $Sub(X)$, and let $\Lambda\subset T^*X$
be an arbitrary conical Lagrangian.
There is a fringed set $R\subset \R^2$ such that for all $(\eta,\kappa)\in R$,
we have
\begin{enumerate}
\item $\eta$ is not a $\Lambda$-critical value of $m$,
\item $\kappa$ is not a $\Lambda$-critical value of $n$,
\item $(\eta,\kappa)$ is not a $\Lambda$-critical value of $m\times n$.
\end{enumerate}
\end{lem}

\begin{proof}
Critical values form a closed $\CC$-subset, 
and their complement is dense.
\end{proof}

First, fix a Whitney stratification $\CS_0$ of $X$ compatible with $Y_0\subset X$,
and let $\Lambda_{\CS_0}\subset T^*X$ be the conical conormal set associated to $\CS_0$.
Apply the preceding lemma to $\fY_1=(Y_1,m_1,n_1)$ and $\Lambda_{\CS_0}$
to obtain a fringed set $R_1\subset \R^2$.
For any $(\eta_1,\kappa_1)\in R_1$, let $T_1\subset X$ be the open submanifold with corners 
$$
T_1 = \{x\in X | m_1(x) > \eta_1, n_1(x)<\kappa_1\}.
$$
We think of $T_1$ as a tube around $Y_1$. We refer to the codimension one
boundary piece 
$$
E_1 = \{x\in X | m_1(x) = \eta_1, n_1(x)<\kappa_1\}
$$
as the end of $T_1$, and the the codimension one
boundary piece 
$$
S_1=\{x\in X | m_1(x) > \eta_1, n_1(x)=\kappa_1\}.
$$
as the side of $T_1$.

Next,
for any $(\eta_1,\kappa_1)\in R_1$, 
fix the Whitney stratification $\CS_{(\eta_1,\kappa_1)}$ of $X$ given by $T_1$, the codimension one
pieces of its boundary, the corners of its boundary,
and the complement of its closure.
Now let $\Lambda_{\CS_{(\eta_1,\kappa_1)}}$ be 
the conical conormal set associated to $\CS_{(\eta_1,\kappa_1)}$.
Apply the lemma to $\fY_0=(Y_0,m_0,n_0)$ and $\Lambda_{\CS_1}$
to obtain a fringed set $R_0\subset \R^2$.
For any $(\eta_0,\kappa_0)\in R_0$, let $T_0\subset X$ be the open tube
$$
T_0 = \{x\in X | m_0(x) > \eta_0, n_0(x)<\kappa_0\},
$$
with end
$$
E_0 = \{x\in X | m_0(x) = \eta_0, n_0(x)<\kappa_0\},
$$
and side
$$
S_0=\{x\in X | m_0(x) > \eta_0, n_0(x)=\kappa_0\}.
$$

We will also need the relative dualizing objects $\omega_{T_0/Y_0},$ $\omega_{T_1/Y_1}$.
To construct these, 
choose retracting 
fibrations of pairs 
$$\ol\pi_0:(T_0\cup E_0,E_0)\to ((T_0\cup E_0)\cap Y_0,E_0\cap Y_0),
\ol\pi_1:(T_1\cup E_1,E_1)\to ((T_1\cup E_1)\cap Y_1,E_1\cap Y_1),
$$
 consider the restrictions $\pi_0 =\ol\pi_0|_{T_0}$, $\pi_1 =\ol\pi_1|_{T_1}$, and
define 
$$\omega_{T_0/Y_0} = \pi_0^! \C_{Y_0}, 
\quad
\omega_{T_1/Y_1}=\pi_1^! \C_{Y_1}.
$$
Concretely,  $\omega_{T_0/Y_0},$ $\omega_{T_1/Y_1}$ are canonically isomorphic
to the local systems (placed in degrees $-\codim Y_0$, $-\codim Y_1$)
on $T_0,$ $T_1$ of relative orientations along the fibers of $\pi_0$, $\pi_1$ respectively.

Finally, we define the morphisms in the dg category $Sub(X)$ to be given by the
relative de Rham complex
$$
\hom_{Sub(X)}(\fY_0,\fY_1) = (\Omega(T_0\cap T_1, E_0 \cup S_1;
\omega^{-1}_{T_0/Y_0}\otimes \omega_{T_1/Y_1}),d).
$$
Given a finite collection of objects of $Sub(X)$,
we may generalize the above perturbation procedure in order to define
the composition of morphisms as the wedge product of differential forms.

\begin{prop} 
For submanifolds $i_0:Y_0\hookrightarrow X$, $i_1:Y_1\hookrightarrow X$, we have
a canonical quasi-isomorphism 
$$
\hom_{Sh(X)}(i_{0*}\C_{Y_0}, i_{1*}\C_{Y_1})\simeq 
(\Omega(T_0\cap T_1, E_0 \cup S_1; 
\omega^{-1}_{T_0/Y_0}\otimes \omega_{T_1/Y_1}),d).
$$
The composition of morphisms coincides with the wedge product of differential forms.
\end{prop}

\begin{proof}
The proof is similar to the proofs of Sections~\ref{open} and~\ref{smooth}.

Consider the inclusions
$$
T_0\stackrel{s_0}{\hookrightarrow} 
T_0\cup S_0\stackrel{e_0}{\hookrightarrow} 
X
\qquad
T_1\stackrel{s_1}{\hookrightarrow} 
T_1\cup S_1\stackrel{e_1}{\hookrightarrow} 
X
$$
Then by de Rham's theorem, we have a quasi-isomorphism
$$
\hom_{Sh(X)}(e_{0*}s_{0!}\omega_{T_0/Y_0},
e_{1*}s_{1!}\omega_{T_1/Y_1})
\simeq
(\Omega(T_0\cap T_1, E_0 \cup S_1;
\omega^{-1}_{T_0/Y_0}\otimes \omega_{T_1/Y_1}),d).
$$

One may identify the left hand side of this quasi-isomorphism with that of the proposition using standard
identities as in Lemma~\ref{hom calc lemma}, and repeated applications of 
the Thom isotopy lemma as in 
Lemmas~\ref{first app of isotopy} and~\ref{second app of isotopy}.
We leave the details including the last assertion to the interested reader.
\end{proof}

By the preceding proposition, we may define a dg functor 
$$P:Sub(X)\to Sh(X)$$
by sending an object $\fY=(Y,m,n)$ to the standard sheaf $i_*\C_Y$ where
$i:Y\hookrightarrow X$ is the inclusion. The induced dg functor on twisted
complexes $\tw P:\tw Sub(X) \to Sh(X)$ is a quasi-equivalence.


\subsection{Morse theory interpretation}

In Section~\ref{morse theory section}, we introduced
the $A_\infty$-category $Mor(X)$ and showed
it is quasi-equivalent to $Open(X)$.
The results of this section generalize that discussion.

We extend the definition of $Mor(X)$ as follows.
As with $Sub(X)$,
we take the objects of $Mor(X)$ to be triples
$(Y,m,n)$ where $Y\subset X$ is a submanifold, $m:X\to \R_{\geq 0}$ is a defining function
for its boundary $\del Y\subset X$, and $n:X\to \R_{\geq 0}$
is a defining function
for its closure $\ol Y\subset X$. 
To define the complex of morphisms
from an object $\fY_0=(Y_0,m_0,n_0)$
to an object $\fY_1=(Y_1,m_1,n_1)$, we introduce some constructions refining those
of the previous section.

To refine the procedure of the preceding section, 
we first fix $(\eta_1,\ol \kappa_1)\in R_1$, 
and consider the function
$$
f_1 = \log m_1 - \log (\ol \kappa_1 -n_1).
$$
For any positive $\kappa_1 <\ol \kappa_1$, we have the open tube
$$
T_1 = \{x\in X | m_1(x) > \eta_1, n_1(x)<\kappa_1\}.
$$
If $\kappa_1$ is sufficiently close to $\ol \kappa_1$,
then there is a convex open set of Riemannian metrics on $X$ for which
the gradient $\nabla f_1$ is inward pointing along the end $E_1$ and outward
pointing along the side $S_1$.

We next proceed similarly and 
choose $(\eta_0,\ol \kappa_0)\in R_0$, 
and consider the function
$$
f_0 = \log m_0 - \log (\ol \kappa_0 -n_0).
$$
For any positive $\kappa_0 <\ol \kappa_0$, we have the open tube
$$
T_0 = \{x\in X | m_0(x) > \eta_0, n_0(x)<\kappa_0\}.
$$
For sufficiently small positive $\eta_0$, and $\kappa_0$ sufficiently close to $\ol \kappa_0$,
there is a convex open set of Riemannian metrics on $X$ for which
the gradient $\nabla f_1 -\nabla f_0$ is inward pointing along the end $E_0\cap T_1 $ and outward
pointing along the side $S_0\cap T_1$.
To insure analogous but opposite behavior along the end $T_0\cap E_1$ and the side $T_0\cap S_1$,
we proceed as follows. 
By moving $f_0$ to a new function $\wt f_0$, 
we may arrange so that 
 there is an open convex set of Riemannian metrics on $X$ for which
the gradient $\nabla f_1 -\nabla \wt f_0$ is outward pointing along the end $T_0\cap E_1$ and inward
pointing along the side $T_0\cap S_1$. Furthermore, 
for an open convex set of Riemannian metrics on $X$, we continue to have that 
$\nabla f_1 -\nabla\wt f_0$ is inward pointing along the end $E_0\cap T_1 $ and outward
pointing along the side $S_0\cap T_1$.

Finally, we choose small perturbations of our functions and metric, and define the morphisms 
of $Mor(X)$ to be the Morse complex
$$
\hom_{Mor(X)}(\fY_0,\fY_1) = 
(\bigoplus_{p\in Cr(T_0\cap T_1, f_1 -\wt f_0)} 
\intHom(\omega_{T_0/Y_0}|_p, \omega_{T_1/Y_1}|_p), \ m^1_{Mor(X)}).
$$
The verification that this is a well-defined complex is similar to Lemma~\ref{morsecomplex}.
As usual, to define the higher compositions, one generalizes the above procedure
sequentially for a finite collection of objects. The details of this are no more complicated than in
other contexts considered earlier. 
Similarly, the fact that we obtain an $A_\infty$-category
follows from homological perturbation theory
along the same lines as the arguments of Section~\ref{homol pert theory applied}.
In addition,
as in Section~\ref{homol pert theory applied},
homological perturbation theory also provides
an $A_\infty$-quasi-equivalence 
$$
M:Sub(X)\to Mor(X).
$$


\subsection{Identification with standard branes}

In Sections~\ref{perturbation section} and ~\ref{relation to morse theory section}, 
we explained how to calculate morphisms in the Fukaya category
among standard branes associated to open submanifolds. In this section, we adapt
that discussion to the case of standard branes associated to arbitrary submanifolds.

Recall that given a defining function $m:X\to \R_{\geq 0}$ for the boundary $\del Y\subset X$,
we define the standard Lagrangian $L_{Y,m}\subset T^*X$ to be the fiberwise sum
$$
L_{Y,m} =T^*_Y X + \Gamma_{d\log m}
$$
where $T^*X_Y\subset T^*X$ is the conormal bundle to $Y$, and $\Gamma_{d\log m}\subset T^*X$
is the graph of the differential of $\log m$.
It comes equipped with a canonical flat bundle and brane structure and thus may be considered 
as an object of $Fuk(T^*X)$.


\subsubsection{Perturbations}

We first explain the necessary modifications to the perturbation procedure of 
Section~\ref{perturbation section}.
Recall that our perturbations were made up of two steps: a variable dilation followed by
a separation at infinity. In our current setting, this may not be enough to guarantee
that the height of our Lagrangians will be small enough along certain critical contours. 
Namely, it may not hold that the pairwise differences of our Lagrangians provide
vector fields with prescribed inward and outward behavior along the codimension one
boundary components
of the intersections of certain open tubes.
Thus we will add a third independent step: a final variable dilation.
Before explaining this, we first comment on the only substantive change in the first two steps.

The first step involving a variable dilation remains the same. But we will change the second step
as follows. Rather than using it solely to separate Lagrangians near infinity, we will also use it to 
tilt our standard branes so that they become very close to being graphs over open subsets.
We will explain this in the setting of a single standard brane $L_{Y,m}$ and leave it to the reader
to repeat the arguments of Section~\ref{perturbation section} using this version of the separation step.

We reinterpret the separation flow as follows. Fix large positive numbers 
$k<\ell \in(0,\infty)$, and choose a decreasing function
$b_{k,\ell}:\R\to\R$ satisfying the following
$$
b_{k,\ell}(r) =
\left \{
\begin{array}{cl}
1 & \mbox{ for $r\leq k$,} \\
0  & \mbox{ for $r\geq \ell$.} \\
\end{array}
\right.
$$
It will be convenient to consider the translation of the cotangent bundle $T^*X$ where the
zero section is given by the scaled Lagrangian
$$
Z= b_{k,\ell}(|\xi|) \cdot \Gamma_{d\log m}
$$ 
Note that though $\Gamma_{d\log m}$ is a singular graph, the function $b_{k,\ell}(|\xi|)$
is constructed to be zero near the singularities. Thus the section $Z$ is a well-defined graph over all of $X$. 
Observe that fiberwise addition by $Z$ is a symplectomorphism.

Fix $\ol\kappa>0$, and consider the function $h:\R\to\R$ defined by
$$
h(r)= 
\int_0^r \frac{-1+\sqrt{1 + 4r^2\ol\kappa}}{2r}
$$
for $r\not = 0$, and $h(0)=0$.
Note that $h$ is differentiable. Consider the Hamiltonian
$$
H(x,\xi)= h(|\xi - Z|)
$$
where we take the length of the fiberwise difference.
The associated vector field $v_H$ is
in the direction $v_\theta$ with length $h'(|\xi- Z|)$. 
Thus the flow $\varphi_{H,t}$ associated to $v_H$ is nothing more than a rescaled version 
of the geodesic flow.

The point of choosing $h$ as we have is the following.
Applying $\varphi_{H,t}$ to the standard Lagrangian 
$\Gamma_{d\log m}$ leads to a perturbation with similar characteristics
as that considered in Section~\ref{separation}.
Applying  $\varphi_{H,t}$ to the standard Lagrangian 
$L_{Y,m}$ produces a perturbed Lagrangian which 
is a graph over an open tube $T$ around $Y$. To be more precise,
recall from the preceding sections the construction of the tube $T$ associated to an object $\fY=(Y,m,n)$.
Namely, we fix $(\eta,\ol \kappa)$ in the fringed set $R$, and then for any positive $\kappa<\ol\kappa$,
we form the open tube
$$
T = \{x\in X | m(x) > \eta, n(x)<\kappa\}.
$$
Now assume within the region
$m\geq\eta$ the defining function $n$ is equal to half the squared-distance from $Y$ . 
(Note that this is not a significant constraint 
since we may choose $m$ and $\eta$ independently beforehand.)
Then over the tube $T$, the unit time perturbation $\varphi_{H,1}(L_{Y,m})$ will be the graph
of the differential of the function
$$
f = \log m - \log (\ol \kappa -n).
$$
In the next section, we will use this compatibility with the previously defined
Morse category $Mor(X)$ in order to see that calculations in $Mor(X)$ agree
with those in the Fukaya category $Fuk(T^*X)$.
Note that though the above perturbation is only asymptotically
normalized geodesic flow, it may easily be modified to be precisely normalized
geodesic flow near infinity without changing any of its other properties.

Finally, we add a third step to our perturbation procedure to ensure 
that the height of our Lagrangians will be small enough along certain critical contours. 
Namely, we must further move our perturbed Lagrangian to
guarantee that in analogy with Proposition~\ref{pertout},
when we consider multiple standard Lagrangians, their pairwise differences
will have the correct behavior (inward and outward pointing)
along the boundaries of the intersections of the corresponding tubes.


\subsubsection{Relation to Morse theory}

The arguments of Section~\ref{relation to morse theory section} 
extend directly to this setting. 
To simplify things, 
we may work with objects $\fY=(Y,m,n)$ 
such that $n$ is equal to half the squared-distance from 
the closure $\ol Y\subset X$ away from the boundary $\del Y\subset \ol Y$. 
More precisely, we may assume $n$ is equal to half the squared-distance from 
$Y$ within the region
$m\geq\eta$, where $(\eta,\kappa)$ is in the fringed set $R$ for small enough $\kappa$.
Then as we have seen, 
our Morse perturbations and Fukaya
perturbations are compatible.
To control the possible holomorphic
polygons, we proceed similarly as in Section~\ref{relation to morse theory section}.
The only amendment is that here for each object $\fY=(Y,m,n)$ 
when we show that polygons do not escape as in Lemma~\ref{second area bound},
we begin the arguments with a region $S$
which is a small neighborhood of the hypersurface with corners
$$
X_{m =\eta', n\leq \kappa'} \cup X_{m\geq  \eta', n = \kappa'}.
$$
Then as before, 
applying the theorem of Fukaya and Oh provides the desired identification
of moduli spaces. 
In conclusion, we obtain an $A_\infty$-quasi-embedding
$Mor(X)\hookrightarrow Fuk(T^*X)$ extending that of Section~\ref{relation to morse theory section},
and that takes the standard object $\fY=(Y,m,n)$ to the standard brane $L_{Y,m}$.


\subsection{Other objects}

We informally mention here another class of objects of $Sh(X)$ which also go to Lagrangians
under our quasi-embedding: the so-called {tilting perverse sheaves}. These
may be thought of as extensions of flat vector bundles on submanifolds
with boundary conditions somewhere between the standard and costandard extensions.
While the intersection cohomology or intermediate extension is cohomologically
between the standard and costandard extensions, the tilting extension (if it exists)
is geometrically between the two.
To understand this, 
note that one can view the standard Lagrangian $L_{Y,m}$ as giving a vector field on $Y$
which is everywhere inward pointing along $\del Y$,
and similarly,
the costandard Lagrangian $-L_{Y,m}$ as giving a vector field 
which is everywhere outward pointing.
The Lagrangians associated to tilting perverse sheaves give vector fields
which are sometimes inward and sometimes outward pointing over prescribed parts of the boundary.

Rather than further developing this picture here, we content ourselves with giving
an example and picking up the discussion elsewhere. Consider the complex
line $\C\simeq\R^2$ with coordinate $z$, and let $i:U\hookrightarrow X$ be the open 
subset $U=\{z\in \C | z\not = 0\}$. 
Given the defining function $m(z)=|z|^2/2$ for 
the point $0\in\C$,
the standard Lagrangian corresponding to the standard sheaf
$i_*\C_U$ is given by the graph of the real part of $dz/z$.
Similarly, the costandard Lagrangian corresponding to the costandard sheaf
$i_!\C_U$ is given by the graph of the real part of $-dz/z$. The graph
of the real part of $dz/z^n$, for $n\geq 2$, is a Lagrangian corresponding to 
a tilting perverse sheaf. In particular, for $n=2$, it corresponds
to the indecomposable tilting extension of the constant sheaf on~$U$.


\vskip 0.2in
\noindent
{\scriptsize
{\bf David Nadler,} Department of Mathematics, Northwestern University,
2033 Sheridan Road, Evanston, IL  60208 (nadler@math.northwestern.edu)\\
{\bf Eric Zaslow,} Department of Mathematics, Northwestern University,
2033 Sheridan Road, Evanston, IL  60208 (zaslow@math.northwestern.edu)
}


\end{document}